\documentclass{article}
\usepackage{amsmath}
\usepackage{amsfonts}
\usepackage{amssymb}
\usepackage{pstricks,pst-node,amssymb,amsmath,graphics,latexsym,tabularx,shapepar}
\usepackage[all,cmtip,2cell,line]{xypic}
\usepackage{graphicx}%
\newtheorem{theorem}{Theorem}

\newtheorem{conclusion}[theorem]{Conclusion}

\newtheorem{corollary}[theorem]{Corollary}

\newtheorem{definition}[theorem]{Definition}
\newtheorem{example}[theorem]{Example}

\newtheorem{lemma}[theorem]{Lemma}

\newtheorem{proposition}[theorem]{Proposition}
\newtheorem{remark}[theorem]{Remark}

\newtheorem{keywords}[theorem]{Keywords}

\newenvironment{proof}[1][Proof]{\noindent\textbf{#1.} }{\ \rule{0.5em}{0.5em}}
\newcommand{\bpartial}{\mathop{\partial\kern -4pt\raisebox{.8pt}{$|$}}}
\newcommand{\bra}{\mathopen{[\kern-1.6pt[}}
\newcommand{\ket}{\mathclose{]\kern-1.5pt]}}
\newcommand{\bbra}{\mathopen{[\kern-2.2pt[\kern-2.3pt[}}
\newcommand{\bket}{\mathclose{]\kern-2.1pt]\kern-2.3pt]}}

\makeindex

\begin{document}

\title{ Optimal control problems on the co-adjoint Lie groupoids }
\author{ Ghorbanali Haghighatdoost \\  Department of Mathematics, Azarbaijan Shahid Madani University,\\Tabriz, Iran \\ e-mail: gorbanali@azaruniv.ac.ir  \\[2ex]
           }

\maketitle

\abstract{In this work we study the invariant optimal control problem on Lie groupoids. We show that any invariant optimal control problem on a Lie groupoid reduces to its co-adjoint Lie algebroid. In the final section of the paper, we present  an illustrative example.}




\maketitle






\noindent
\section{Introduction}

 Control systems and optimal control theory are the main branches of most area of the sciences. Especially in mathematics, they have taken a major place. The interplay between optimal control theory, geometry
and analytical mechanics is interest of many authors. Geometric investigations of the optimal control problems are done by many control theorists. Ever since the optimal control  emerged in Mathematics, it has strongly influenced geometry. In particular, it played a key role in the birth of differential geometry and defining the 'straight line' or the concept of 'geodesic' by means of an extermal curve of the problem. But more recently, modern control theory has been heavily influenced by geometry and Hamiltonian mechanics(see \cite{suss, Jur, BW} for more details).

Most works concern invariant control systems. It means that the system has a symmetry, i.e. some group is the symmetry of the system. Many authors studied optimal control problems defined on some Lie groups. This formalism has been taken in V. Jurdjevic's works in \cite{Jur, Jur2}. In this case the theory of Hamiltonian systems on Lie groups is based on a particular realization of the cotangent bundle of a Lie group $G$ as the product of $G$ and the dual of its Lie algebra $\mathfrak{g}$, i.e. $T^*G$ realized as $G\times \mathfrak{g}^*.$

After that  is considered some Hamiltonian function for every invariant vector field on the realization of the co-tangent bundle of the Lie group and are formed the reduced equations which are defined in terms of Lie-Poisson structure on the dual of the Lie algebra of the Lie group. In the following, by a suitable symplectic form is defined correspondence Hamiltonian vector field and Hamiltonian equation and proved that the extermal curves of the optimal control problem are among of the integral curves of this Hamiltonian system. 

The generalization of the Jurdjevic's idea is done by E. Martinez in \cite{Mar}, where control systems and optimal control problems on a Lie algebroid are studied. By considering the prolongation of the dual bundle of the Lie algebroid with respect to the Lie algebroid itself is introduced a Hamiltonian function on the direct product of the dual bundle of the Lie algebroid and a bundle so-called control space over a base manifold of the Lie algebroid.

The Hamiltonian section associated with the Hamiltonian function is constructed and is shown that integral curves of the Hamiltonian vector field, which is given by the projection of the Hamiltonian section under the anchor of the prolongation, are the critical trajectories. Then, E. Martinez  showed that the solutions of the optimal control problem are described by trajectories of the Hamiltonian system.

It is worth noting that in \cite{Jur, Jur2} have the optimal control problems on the Lie group and on the co-adjoint orbits of the Lie group. The idea is continued by realization of the co-tangent bundle of the Lie group and by defining a symplectic form over the co-tangent bundle and constructing the Hamiltonian equations. It is shown that the solutions of these equations are the solutions of the optimal control problem.

Unfortunately, for the Lie groupoid, the co-tangent bundle does not admit any realization, so we cannot apply the Jurdjevic's method for the Lie groupoid and get extermal solutions. But in \cite{Joz, GJ} it is shown that any control system as well as any optimal control problem on a Lie groupoid reduces to its Lie algebroid.

In this work, at first we study a right-invariant control system as well a right- invariant optimal control problem on an arbitrary Lie groupoid and show that they can be reduced to the co-adjoint Lie groupoid of the Lie groupoid which  has been constructed in \cite{HA1}. Then by using the reduction in the co-adjoint Lie algebroid and by applying the method which is presented in \cite{Mar}, we study control systems and optimal control problems on the co-adjoint Lie algebroid.
We end the paper with an illustrative example. In other words, we consider the trivial Lie groupoid and  an optimal control problem on its co-adjoint Lie algebroid and  show that the optimal control problem can be reduced to the optimal control problem on the co-tangent bundle of  the orbits of the co-adjoint representation of the Lie group. Also, we show that the extermal solutions of the optimal control problem on the co-adjoint Lie algebroid of the trivial Lie groupoid are obtained from the solutions of the corresponding Hamiltonian system on the co-tangent bundle of the co-adjoint orbits of the Lie group.  
The paper is organized as follows. In section 2 we recall some facts about co-adjoint Lie groupoid and co-adjoint Lie algebroid (see \cite{HA1} for more details) and define a control system on  a Lie groupoid. In section 3 we study reduction in optimal control problem on a Lie groupoid to its co-adjoint Lie groupoid. In section 4 we study the reduction of  the optimal control problem on a Lie groupoid to its Lie algebroid. In section 5 we study the optimal control problem on the co-adjoint Lie algebroid of a regular Lie groupoid. The paper, in section 6, is finished by an illustrative example.

\section{ Main Concepts }
\subsection{ Lie groupoid and Lie algebroid }
It is well-known that a  groupoid  which is  denoted by $G \rightrightarrows M$, consists of two sets $G$ and $M$ together with structural mappings $\alpha, \beta, 1, \iota$ and $m,$  where source mapping $\alpha: G \rightarrow M,$ target mapping $\beta:G \longrightarrow M,$ unit mapping $1: M \longrightarrow G,$ inverse mapping $\iota : G \longrightarrow G $ and multiplication mapping $ m: G_{2} \longrightarrow G $ where $ G_{2} = \lbrace (g,h) \in G \times G ~\vert ~ \alpha (g) = \beta (h) \rbrace $ is subset of $G \times G.$\\

A Lie groupoid is a groupoid $G \rightrightarrows M$ for which $G$ and $M$ are smooth manifolds, $\alpha, \beta, 1, \iota$ and $m,$ are differentiable mappings and besides of $\alpha,\beta$ are differentiable submersions.

The right translations on a Lie groupoid $G$ over $M$, $R_g: G_{\beta(g)}=\alpha^{-1}(\beta(g))\to G_{\alpha(g)}=\alpha^{-1}(\alpha(g))$ are diffeomorphisms of the $\alpha$-fibers only and not of the whole groupoid. 

A smooth mapping $X: G\to TG$ is called a vector field on $G$, i.e. for every $g\in G$, $X(g)\in T_gG$, where $T_gG$ is the tangent space  to  $G$ at $g\in G$.

According to the above, to talk about right-invariant vector fields on $G$, we have to restrict attention to those vector fields which are tangent to the $\alpha$-fibers. In other words, we take the elements of the sections of the sub-bundle $T^\alpha G$ of $TG$ defined as $$T^\alpha G=\mathrm{Ker}(d\alpha)\subset TG.$$

A Lie algebroid $A$ over a manifold $M$ is a vector bundle  $\tau : A \longrightarrow M$ with the following items:
\begin{enumerate}
\item A Lie bracket $[\vert ~,~ \vert] $ on the space of smooth sections $ \Gamma( \tau),$
$$[\vert ~,~ \vert]: \Gamma (\tau) \times \Gamma (\tau) \longrightarrow \Gamma (\tau), \quad (X , Y ) \longmapsto [\vert X , Y \vert].$$
\item A morphism of vector bundles $\rho : A \longrightarrow TM,$ called the anchor map, where $TM$ is  the tangent bundle of $M,$ such that the anchor and the bracket satisfy the following Leibniz rule:
$$[\vert X , f Y \vert] = f [\vert X , Y \vert] + \rho (X) (f) Y,$$
\end{enumerate}
where $X , Y \in \Gamma (\tau)$, $~f \in C^{\infty} (M)$ and $\rho (X) f$ is the derivative of $f$  along the vector field $\rho (X).$

Given a Lie groupoid $G$ over $M$, we define the vector bundle $A=Lie(G)=AG$, whose fiber at $x\in M$ coincides with the tangent space at the unit $1_x$ of the $\alpha$- fiber at $x$. In other words, $AG:=(T^\alpha G)_M.$

It is easy to see that every fiber of the sub-bundle $T^\alpha G$  at an arrow  $h: y\to z$ of $G$ is $T^{\alpha}_h G= T_hG(y,-)$,  where  $G(y,-) =\alpha^{-1}(y) = G_{\alpha(h)}$. Consider the right translation $R_g:G(y,- ) \to G(x, -), g' \to g'g$.
The differential of the right translation by $g$ induces a map $ dR_g: T^\alpha_h G \to T^\alpha_{hg} G$.
\begin{definition}
Vector field $X$ on $G$ is called vertical if it is vertical with respect to $\alpha$, that is, $X_g\in T_gG_{\alpha(g)}$, for all $g\in G$. We call  $X$  right- invariant on $G$ if it is vertical and $X_{gh}=dR_g(X_h)$, for all $(h,g)\in G_{(2)}$. 
\end{definition}
It is easy to show that  $\Gamma(AG)$ - the space of sections of vector bundle $AG$ can be identified the space of right-invariant vector field on $G$. We denote the space of right-invariant vector field on $G$ by 
$$\chi^\alpha_{\mathrm{inv}}(G)=\{X\in \Gamma(T^\alpha G): X_{hg}=dR_g(X_h), \;(h, g)\in G_{(2)}\}.$$
From above, we have the space of sections  $\Gamma (AG)$ which is isomorphic to the space of right-invariant vector fields on $G$, $\chi^\alpha_{\mathrm{inv}}(G)$. 
On the other hand, the space $\chi^\alpha_{\mathrm{inv}}(G)$ is a Lie sub-algebra of the Lie algebra $\chi(G)$ of  vector fields on $G$ with respect to the usual Lie bracket of vector fields. 
Also, the push-forward of vector field on the $\alpha-$fibers along $R_g,$ preserves brackets.
 So we obtain a new bracket on $\Gamma(AG)$ which is uniquely determined. The Lie bracket on $AG$ is the Lie bracket on  $\Gamma(AG)$ obtained from the Lie bracket on $\chi^\alpha_{\mathrm{inv}}(G)$. 
 The anchor of $AG$ is the differential of the target mapping $\beta$, i.e. $\rho=(T\beta)_{AG}: AG \to TM$. 
 As a result, we obtain that $AG$ is a Lie algebroid associated to the Lie groupoid $G$.

\subsection{Control system}
A curve $g: I\to G$ is called $\alpha$-curve if $g$ be smooth and $\frac{dg}{dt}\in T^\alpha G$, i.e. for all $t\in I,$ $g(t)$ belongs to the $\alpha$-fibers of $G$.

A control system on a Lie groupoid $G$ is an ODE of the form
\begin{equation}
  \frac{dg}{dt}= F(g, u).
  \label{cs}
\end{equation}
Here $F: G \times \mathbb R^m \to  T^\alpha G$ is a smooth mapping and  $u(t)=(u_1(t), u_2(t),..., u_n(t))$ taking values in a subset $U$ of
 $\mathcal {R}^m$ is called a control function.

The control system (\ref{cs}) is said right-invariant if for all $(h, g) \in G_{(2)}$ and $u \in U$, we have $$dR_g(F(h, u))=F(R_g(h), u).$$
It is easy to show that every right-invariant vector field on a Lie groupoid $G$ over the manifold $M$ is determined by its value at the points in $M$. Also, each right-invariant system is uniquely defined by its values at the points in $M$. In other words, 
$$ dR_g(F(1_{\beta(g)}, u))=F(R_g(1_{\beta(g)}), u)= F(g, u).$$

By a control system on a Lie groupoid $G$ over $M$ we mean a system of differential equations (\ref{cs}) where $g$ is said to be the state point and $u$ are the control coordinates. Solutions of the system (\ref{cs}) are said to be trajectories of the system. \\

An optimal control problem consists in finding the trajectories of the control system which connect some predetermined states and minimize the integral of some function $f$  so-called cost function depending on state and control coordinates.\\

A function $f: G\times \mathbb{R}^m\to \mathbb R$ is called right-invariant if $$f(R_g(h), u)=f(h, u),$$   for all $(h, g) \in G_{(2)}$ and $u \in \mathbb R^m$.

Every right-invariant function on $G\times \mathbb R^m$ is uniquely determined by its values at the points in $M$. In other words, $f(g,  u)=f(R_g(1_{\beta(g)}),  u)$.

An optimal problem on $G$ is said to be right-invariant if both cost function $f$ and the control system are right-invariant.

\subsection{Co-adjoint Lie groupoid and Co-adjoint Lie algebroid}
A Lie groupoid $ G \rightrightarrows M $ together with structural mappings $\alpha, \beta, 1, \iota$ and $m,$ which  are defined above, is called regular if the mapping $(\beta, \alpha)$ is mapping with constant rank. In \cite{HA1} we associated  to every regular Lie groupoid $G$  over $M$ a Lie groupoid so-called co-adjoint Lie groupoid. The adjoint and co-adjoint action of $G$ on the isotropy Lie groupoid $I_{G}$ of $G$, are defined as follows:

The isotropy group of $ G \rightrightarrows M $ at $p\in M$ is  $ I_{p} = \alpha ^{-1} (p) \cap \beta ^{-1} (p)$. It is well-known that $I_{p}$ is a Lie group where its composition law and inverse map are the restrictions of multiplication map $m$ and inverse map $\iota $ to $I_{p},$ respectively. The union of all isotropy groups $I_{p}$ when $p$ rounds over in $M$ construct a groupoid over $M$, i.e. $I_{G} = ( \cup I_{p} )_{p \in M}$ is a groupoid over $M$. 
Note that the isotropy groupoid of Lie groupoid is not a smooth manifold in general.
If we consider $G$ being a regular Lie groupoid, then its associated isotropy groupoid is Lie groupoid. We denote the associated isotropy Lie groupoid to Lie groupoid $G \rightrightarrows M$ by $I_{G}$ and the Lie algebroid associated to isotropy Lie groupoid by $A I_{G}$ and call it isotropy Lie algebroid.\\

\begin{definition}
\label{groupoid action}
A smooth left action of Lie groupoid G on smooth map  $J : N \longrightarrow M$ is a smooth map $\theta : G_{~ \alpha} \times_{~J} N \longrightarrow N $ which satisfies the following properties:
\begin{enumerate}
\item For every $(g,n) \in G_{~ \alpha} \times_{~J} N,~~~~ J (g.n) = \beta (g),$ 
\item For every $n \in N,~~~~ 1_{J(n)} .n = n,$
\item For every $(g , g^{\prime}) \in G_{2} $ and $n \in J^{-1} (\alpha (g^{\prime})),~~~~ g . ( g^{\prime} . n ) = (g g^{\prime})  . n$
\end{enumerate}

(where $ g.n := \theta ( g , n )$ and $ \theta (g) (n):= \theta (g , n )$).\\

\end{definition}

Similar to the action of the Lie groupoid on smooth mapping (definition \ref{groupoid action}), the definition of action of a Lie algebroid on a smooth mapping will be as follows:
\begin{definition}
An action of a Lie algebroid $ ( A, M, \pi, \rho , [\vert ~,~ \vert ] )$ on map $J : N \longrightarrow M$ is a map\\
$ \theta : \Gamma (A) \longrightarrow \mathfrak{X} (N) $ which for all $ f \in C^{\infty} (M)$ and $ X, Y \in \Gamma^{\infty} (A),$ satisfies in the following properties:
 
\begin{enumerate}
\item $ \theta ( X + Y ) = \theta (X) + \theta (Y)$
\item $ \theta ( f X ) = J^{\ast} f \theta (X) $
\item $ \theta  ( [\vert X , Y \vert] ) = [ \theta (X) , \theta (Y) ]  $
\item $T J ( \theta (X) ) = \rho (X)$
\end{enumerate}
where $J^{\ast} : C^{\infty} (M) \longrightarrow C^{\infty} (N) $ such that $J^{\ast} f = f \circ J \in C^{\infty} (N)$ is pullback of $f$ by $J$.
\end{definition}

\begin{remark}
\label{action remark}
Let $ \theta$ be the action of a Lie groupoid $G$ on smooth map  $J: N \longrightarrow M$ which was introduced into definition \ref{groupoid action}.
As mentioned in \cite{Bos},  every action of a Lie groupoid $G$ on $ J : N \longrightarrow M $ induces an action $\theta^{\prime}$ of a Lie algebroid $ A(G) $ on $ J : N \longrightarrow M $ as follows:
$$ \theta ^{\prime} (X) (n) := \dfrac{d}{dt} \vert _{t=0} ~Exp ( t X )_{J(n)} . n.$$
\end{remark}

\begin{definition}
Let $G$ be a Lie groupoid over $M.$ Then $G$ is regular Lie groupoid if the anchor $( \beta , \alpha ) : G \longrightarrow M;~ g \longmapsto \big( (\beta (g) , \alpha(g) \big)$ is a mapping of constant rank.

\end{definition}

Consider the regular Lie groupoid $ G \rightrightarrows M $ and its associated isotropy Lie groupoid $I_{G}.$  $G$ acts smoothly from the left on $ J: I_{G} \longrightarrow M $ by conjugation, means  $ C:G \times I_{G} \longrightarrow I_{G},$ $~~ C(g) (g^{\prime} ):= g g^{\prime} g^{-1}$ is an action of $G$ on $I_{G}$ which we call it  conjugation action.\\

On the other hand, the  conjugation action induces an action of a Lie groupoid $G$ on $ A I_{G} \longrightarrow M $. We call this action \textbf{adjoint action} of $G$ on $A I_{G}$ which can be defined as follows:
$$\mathrm{Ad}: G \times A I_{G} \longrightarrow AI_{G}, $$
$$\mathrm{Ad}_{g} X :=(\frac{d}{dt})_{t=0}~ C(g) \mathrm{Exp} (tX),$$
where  $p \in M$,   $g \in G_{p} = \alpha^{-1} (p) $ (the $\alpha$-fibers over $p$)   and   $X \in (A I_{G})_p.$\\

The action $\mathrm{Ad}$ induces an adjoint action of $AG $ on $ A I_{G}  \longrightarrow M $ as follows:
$$\mathrm{ad}: AG \times AI_G \to AI_G,$$
$$ \mathrm{ad}_XY=\mathrm{ad} (X)  (Y) := \frac{d}{dt} \Big\vert_{t=0}~ \mathrm{Ad} \left(\mathrm{Exp} (tX) \right) Y, $$
where $ X \in (AG)_{p} ,~ Y \in (A I_{G})_{p}$ and $ p \in M.$ \\

One can easily prove that for every $ X \in \Gamma(AG) $ and $Y \in  \Gamma(A I_{G} )$
$$ \mathrm{ad}_X (Y) = [\vert X , Y \vert].$$ 

Another action of $G$  on dual bundle $A^{\ast} I_{G} $ which is called \textbf{co-adjoint action} of $G,$ is defined as follows:
$$ \mathrm{Ad}^{\ast} : G \times A^{\ast} I_{G} \longrightarrow A^{\ast} I_{G},$$
$$ \mathrm{Ad}_{g} ^{\ast} \xi (X) := \xi ( \mathrm{Ad}_{g^{-1}} X ).  $$
In other words $$ \langle \mathrm{Ad}_{g} ^{\ast} \xi , X \rangle = \langle \xi , \mathrm{Ad}_{g^{-1}} X\rangle, $$\\
where $ g \in G_{p},~ \xi \in (A^{\ast} I_{G})_p.$\\

Again,  the action $\mathrm{Ad}^{\ast}$ induces so-called co-adjoint action of a Lie algebroid $AG$ on $A^{\ast} I_{G} $ which is defined by:
$$ \mathrm{ad}^{\ast} : AG \times A^{\ast} I_{G}  \longrightarrow A^{\ast} I_{G}, $$
$$ \mathrm{ad}_{X} ^{\ast} \xi (Y) := \xi ( \mathrm{ad}_{-X} (Y) ) = \xi ( [\vert Y,X \vert] ) $$
or $$\langle \mathrm{ad}_{X} ^{\ast} \xi , Y  \rangle = \langle \xi , \mathrm{ad} (- X) Y \rangle,$$
where $ \xi \in (A^{\ast} I_{G})_p.$  For more details, about adjoint and co-adjoint actions see \cite{Bos}. \\

Now we define  the co-adjoint Lie groupoid as follows:
$$O(\xi) = \lbrace \mathrm{Ad}_{g} ^{\ast} \xi ~\vert ~ g \in G \rbrace,$$
where $\xi$ is an element of $(A^{\ast} I_{G})_{p}$. It turns out that for all $\xi$ that the stabilizer $G_\xi=\{g : Ad^\ast_g \xi=\xi  \}$ is a normal Lie subgroupoid of $G$, the co-adjoint orbit  $O(\xi)$ has a natural structure of a Lie groupoid with following structural mappings:
 $\alpha^{\prime},\beta^{\prime} , m^{\prime} , 1^{\prime}$ and $\iota^{\prime}$ which are 
given by 
\begin{enumerate}
\item source mapping: $ \alpha^{\prime}: O(\xi) \longrightarrow M, \quad \mathrm{Ad}_{g} ^{\ast} \xi \longmapsto \alpha (g), $
\item target mapping: $ \beta^{\prime}: O(\xi) \longrightarrow M, \quad\mathrm{Ad}_{g} ^{\ast} \xi \longmapsto \beta (g), $
\item multiplication mapping: $m^{\prime} : ( O(\xi) )_{2} \longrightarrow  O(\xi)$ 
 $$(\mathrm{Ad}_{g} ^{\ast} \xi ,\mathrm{Ad}_{h} ^{\ast} \xi ) \longmapsto \mathrm{Ad}_{m(g,h)} ^{\ast} \xi = \mathrm{Ad}_{gh} ^{\ast} \xi, $$
 As we assumed that the Lie subgroupoid $G_\xi$ is normal, so the multiplication $m^{\prime}$ will be  well-defined.
\item unit mapping: $ 1^{\prime}: M \longrightarrow O(\xi), \quad p \longmapsto \mathrm{Ad}_{1_{p}} ^{\ast} \xi , $
\item  inverse mapping: $ \iota^{\prime} : O(\xi) \longrightarrow O(\xi), \quad \mathrm{Ad}_{g} ^{\ast} \xi \longmapsto \mathrm{Ad}_{g^{-1}} ^{\ast} \xi.$
\end{enumerate}
We call $\alpha^{\prime} , \beta^{\prime} , m^{\prime}, 1^{\prime}$ and $\iota^{\prime},$ source, target, multiplication, unit and inverse mapping, respectively, for Lie groupoid $O(\xi).$

We call this Lie groupoid 
\textbf{co-adjoint Lie groupoid} and denote it by $\mathcal{G}_\xi$. In the latter, if it is not confused, for a fixed selected  $\xi$ we denote the co-adjoint Lie groupoid associated to $\xi$ by $\mathcal{G}$.(For more details, see \cite{HA1}).\\

Also it is shown in \cite{HA1} that the Lie algebroid of the co-adjoint Lie groupoid which we call co-adjoint Lie algebriod is
 $$ A \mathcal{G} = \mathrm{Ker} T\alpha^{\prime} \vert_{\mathrm{Ad}_{1_{p}} ^{\ast} \xi} = T_{\xi} O(\xi) \vert_{\mathrm{Ad}_{1_{p}} ^{\ast} \xi}= \lbrace \mathrm{ad}_{X_p} ^{\ast} \xi ~ \vert ~ X_p \in (AG)_p \rbrace. $$

As we mention above, the space of sections of vector bundle
 $A\mathcal{G}$
  can identify the space of right-invariant vector fields on
  $\mathcal{G}$, 
  on the other hand, each right-invariant vector field on the co-adjoint Lie groupoid is determined by a right-invariant vector field on the Lie groupoid $G$. 
  In other words, letting $\mathcal{G}= \lbrace \mathrm{Ad}_{g} ^{\ast} \xi ~\vert ~ g \in G \rbrace$, it is easy to check that, $T^{\alpha'} \mathcal{G}== \lbrace \mathrm{ad}_{X} ^{\ast} \xi ~\vert ~ X \in T^\alpha G \rbrace$. 
  Also by definition of $\alpha'$  for every $y \in M$, $\mathcal G(y, -)=G(y, -)$, therefore,  for  $\tilde{g}= \mathrm{Ad}_{g}^{\ast}\xi \in \mathcal{G}$. We have $R_{\tilde{g}}: \mathcal{G}_{\beta'(\tilde{g})} \to \mathcal{G}_{\alpha'(\tilde{g})}$, $\tilde{h} \to \tilde{h}.\tilde{g}$ or 
$$R_{\tilde{g}}(\tilde{h})=R_{\mathrm{Ad}_{g}^{\ast}\xi}(\mathrm{Ad}_{h}^{\ast}\xi)=\mathrm{Ad}_{g}^{\ast}\xi.\mathrm{Ad}_{h}^{\ast}\xi,$$
where $(h, g)\in G_{(2)}$ and by multiplication in $\mathcal{G}$ we have $\mathrm{Ad}_{g}^{\ast}\xi.\mathrm{Ad}_{h}^{\ast}\xi=\mathrm{Ad}_{hg}^{\ast}\xi= \mathrm{Ad}_{R_{g}(h)}^{\ast}\xi.$ As a result we obtain $$R_{\tilde{g}}(\tilde{h})=\mathrm{Ad}_{R_{g}(h)}^{\ast}\xi,$$ 
where $R_g: G_{\beta(g)}\to G_{\alpha(g)}.$\\
Also, for  $dR_{\tilde{g}}: T_{\tilde{h}}^{\alpha^\prime}\mathcal{G} \to T_{\tilde{h}\tilde{g}}^{\alpha^\prime}\mathcal{G}$,  we have 
\begin{equation}
dR_{\tilde{g}}(\tilde{X}_{\tilde{h}})=\mathrm{ad}_{dR_g (X_h)}^{\ast}\xi,
\label{dr}
\end{equation} 
where $\tilde{X}_{\tilde{h}}=\mathrm{ad}_{X_h}^{\ast}\xi$.
As a result from (\ref{dr}) we obtain the  following lemma:
\begin{lemma}
The vector field $\tilde{X}=\mathrm{ad}_X^{\ast} \xi$ on $\mathcal G$ is right-invariant if only if $X$ is right-invariant on $G$, that is 

$$ dR_{\tilde{g}}(\tilde{X}_{\tilde{h}})=\tilde{X}_{\tilde{h}\tilde{g}} \quad  \text{if \; only \; if}   \quad  dR_{g(X_h)}= X_{hg}.$$
\label{riv}
\end{lemma}
\section{Reduction in optimal control problem}
Consider the invariant control system (\ref{cs}) on a Lie groupoid $G\rightrightarrows M$. In this section we define a control system on the co-adjoint Lie groupoid $\mathcal {G} \rightrightarrows M$.\\
Control system on the co-adjoint Lie groupoid we define as follows:

Let $\tilde{g}=\mathrm{Ad}_g^{\ast}\xi$ be an element of $\mathcal{G}$ and  let $\tilde{g}:I\to \mathcal{G}$ be a curve in $\mathcal{G}$. It is easy to show that $\tilde{g}: I\to \mathcal{G}$ is $\tilde{\alpha}$-curve in $\mathcal{G}$ if  $g: I\to G$ is $\alpha$-curve in $G$.\\
For any  $\tilde{\alpha}$-curve $\tilde{g}$ in $\mathcal{G}$ we consider following system
\begin{equation}
\frac{d\tilde{g}}{dt}= \tilde{F} (\tilde{g}, u),
 \label{cs2}
\end{equation}
where $ \tilde{F} (\tilde{g}, u)=\mathrm{ad}_{F(g,u)}^{\ast}\xi$.
By considering lemma \ref{riv}, we have following main theorem:
\begin{theorem}
The control system (\ref{cs2}) on the co-adjoint Lie groupoid $\mathcal{G}$ is right-invariant if only if the control system (\ref{cs}) is right-invariant on the Lie groupoid $G$.
\label{ri}
\end{theorem}
\proof
Let the system (\ref{cs}) be right- invariant, that is, for all $(g, h) \in G_{(2)}$, i.e. $dR_g \circ F(h, u)=F(R_g(h), u).$

We show that for all $(\tilde{g}, \tilde{h}) \in \mathcal {G}_{(2)}$ and $\tilde{F}(\tilde{h}, u)=\mathrm{ad}_{F(h, u)}\xi$, where $\tilde{g}=\mathrm{Ad}^\ast_g\xi$ and $\tilde{h}=\mathrm{Ad}^\ast_h\xi:$   
$$dR_{\tilde{g}}\circ \tilde{F}(\tilde{h}, u)=\tilde{F}(R_{\tilde{g}}(\tilde{h}), u).$$
We have 
$$dR_{\tilde{g}}\circ \tilde{F}(\tilde{h}, u)=dR_{\tilde{g}}(\mathrm{ad}^\ast_{F(h, u)}\xi)=  \mathrm{ad}^\ast_{dR_g\circ F(h, u)}$$
$$=\mathrm{ad}^\ast_{F(hg, u)}\xi= \tilde{F}(\tilde{hg}, u)=\tilde{F}(\tilde{h}\tilde{g}, u)= \tilde{F}(\mathrm{Ad}^\ast_{hg}\xi, u)$$
$$=\tilde{F}(\mathrm{Ad}^\ast_h\xi.\mathrm{Ad}^\ast_g \xi, u)=\tilde{F}(\tilde{h}\tilde{g}, u)= \tilde{F}(R_{\tilde{g}}(\tilde{h}), u).$$

Conversely, let the system (\ref{cs2})  be right- invariant, i.e. 
$$dR_{\tilde{g}}\circ \tilde{F}(\tilde{h}, u)=\tilde{F}(R_{\tilde{g}}(\tilde{h}), u).$$
At first we prove following relation:
$$dR_{\tilde{g}}(\mathrm{ad}^\ast_X\xi)=\mathrm{ad}^\ast_{dR_{g}(X)}\xi,$$
where $X\in AG$. 
Let $\tilde{X}=\mathrm{ad}^\ast_X\xi \in A\mathcal{G}$, there exists an $\tilde{\alpha}$-curve $\tilde{\gamma}(t) \in \mathcal{G}$, such that $\dot{\tilde{\gamma}}(0)=\tilde{X}$ and $\tilde{\gamma}(0)=\tilde{g}.$ As we describe above 
$dR_{\tilde{g}}(\tilde{X})\in T^{\tilde{\alpha}}_{\tilde{h}\tilde{g}}\mathcal{G}$ , so for all $\tilde{\phi}:\mathcal{G}_{\tilde{\alpha}} \to \mathbb{R}$ we have
$$dR_{\tilde{g}}(\tilde{X})(\tilde{\phi})= \tilde{X}(\tilde{\phi} \circ R_{\tilde{g}})=(\frac{d}{dt})_{t=0}(\tilde{\phi} \circ R_{\tilde{g}}(\tilde{\gamma}(t)))$$
$$=(\frac{d}{dt})_{t=0}\tilde{\phi}(\tilde{\gamma}(t). \tilde{g})= (\frac{d}{dt})_{t=0}\tilde{\phi}(\mathrm{Ad}^\ast_{\gamma(t)}\xi. \mathrm{Ad}^\ast_g \xi),$$
where $\gamma(t)$ is an $\alpha$-curve in $G$ and $g\in G.$ Therefor we have (by multiplication in $\mathcal{G}$)
$$dR_{\tilde{g}}(\tilde{X})(\tilde{\phi})=(\frac{d}{dt})_{t=0}\tilde{\phi}(\mathrm{Ad}^\ast_{\gamma(t).g} \xi)$$
$$= (\frac{d}{dt})_{t=0}\tilde{\phi}(\mathrm{Ad}^\ast_{R_g(\gamma(t))} \xi)= \mathrm{ad}^\ast_{dR_g(X)}\xi(\tilde{\phi}),$$
where $X=\dot{\gamma}(0).$\\
We obtain $$dR_{\tilde{g}}(\tilde{X})=\mathrm{ad}^\ast_{dR_g(X)}\xi. $$
Now if we take as a $\tilde{X}=\tilde{F}(\tilde{h}, u)=\mathrm{ad}^\ast_{F(h, u)}\xi$, we have
$$dR_{\tilde{g}}(\tilde{F}(\tilde{h}, u))=\mathrm{ad}^\ast_{dR_g \circ F(h, u)}\xi, \quad (I)$$
but $$dR_{\tilde{g}}(\tilde{F}(\tilde{h}, u))= \tilde{F}(R_{\tilde{g}}(\tilde{h}), u)= \mathrm{ad}^\ast_{F(R_g(h), u)}\xi. \quad (II)$$
Because $\mathrm{ad}^\ast$ is linear, so from (I) and (II) we get:
$$dR_g \circ F(h, u)= F(R_g(h), u).$$
In other words, the system (\ref{cs}) is right-invariant. So the proof is completed.

 \begin{corollary}
 Every right-invariant control system on any regular Lie groupoid  reduces to its co-adjoint Lie groupoid.
 \label{rcs}
 \end{corollary}

\subsection{Optimal control problem on a Lie groupoid and its reduction to the co-adjoint Lie groupoid}

By definition, an optimal control problem on a Lie groupoid $G\rightrightarrows M$ for the control system (\ref{cs}) is associated with the cost function $f: G\times U \to \mathbb R$ and a controlled pair $(g(t), u(t))$.
 The problem is to find a controlled pair $(g(t), u(t))$ with $t\in [t_0, t_1]$ such that the integral $\int_{t_0}^{t_1}f(g(t), u(t))dt$ is minimal among all controlled pair $(g(t), u(t))$ with some additional conditions that will not be mentioned here.
 Now let us have an optimal control problem on a regular Lie groupoid $G\rightrightarrows M$. This problem is reduced to an optimal control problem on the co-adjoint Lie groupoid $\mathcal{G}\rightrightarrows M$.
 By the corollary \ref{rcs}, every right-invariant control system on a Lie groupoid is reduced to a right-invariant control system on its co-adjoint Lie groupoid. Now if we have a right-invariant cost function $f:G\times U \to \mathbb R$,  we define a new cost function on $\mathcal{G}$ as follows:\\
Consider  function $\tilde{f}: \mathcal{G}\times U \to  \mathbb R$ as $\tilde{f}(\tilde{g}, u)=f(g, u)$ where $\tilde{g}=\mathrm{Ad}_g^{\ast}\xi$. 

If $f$ is right-invariant cost function on $G\times U$ , $\tilde{f}$  is also right-invariant on $\mathcal{G}\times U$ because $$\tilde{f}(R_{\tilde{g}}(\tilde{h}),u)= \tilde{f}(\tilde{h}\tilde{g}, u)=  \tilde{f}(\widetilde{hg}, u)= f(hg, u)=f(R_g(h), u)=f(h, u)=\tilde{f}(\tilde{h}, u),$$
where $\tilde{g}=\mathrm{Ad}_g^{\ast}\xi$ and $\tilde{h}=\mathrm{Ad}_h^{\ast}\xi$.
 
According to above, theorem \ref{ri} and corollary  \ref{rcs}, we have the following proposition.
\begin{proposition}

Any right-invariant optimal control problem on a Lie groupoid $G\rightrightarrows M$  reduces to an optimal control problem on its co-adjoint Lie groupoid $\mathcal{G}\rightrightarrows M.$
\end{proposition}

\section{Optimal control problem on a Lie algebroid}
Some concepts and definitions in this section have been taken from \cite{Joz} and \cite{Mar}.\\
Let $(A, M, \rho, [.,.], \pi)$  be  a Lie algebroid. 


\begin{definition}
A control system on a Lie algebroid $A$ over $M$ with a vector bundle so-called control space $\pi: C \to M$ is a section $f$ of $A$ along $\pi$, i.e. $f: C \to A$. A trajectory of the system $f$ is an integral curve of the vector field $\rho(f)$ along $\pi$, i.e. a trajectory is a solution of following equation
\begin{equation}
\dot{x}(t)= \rho(f(c(t))),
\label{cs4}
\end{equation}
 where $c(t)$ is a curve in $C$ and $x(t)=\pi(c(t))$.
 \end{definition}

Now, if we have an optimal control problem on the Lie algebroid $A$, in other words,  given a cost function $L\in C^\infty (C)$, the problem is to minimize the integral of $L$ over the set of trajectories of the system (\ref{cs4}), i.e. minimize $\int_{t_0}^{t_1}L(c(t))dt$
which $c(t)$ is a trajectory of the system $f$  satisfying some boundary conditions. To solve the problem, a Hamiltonian function $H\in C^\infty(A^\ast \times_M C)$,  $H(\eta, c)=\langle\eta, f(c)\rangle - L(c)$ is used  and its associated  Hamiltonian control system $\sigma_H$ on a subset of $\mathcal{T} A^\ast$-prolongation of $A^\ast$ along $pr_1:A^\ast\times_M C\to A^\ast.$ Because of the Pontryagin Maximum Principle, which is the fundamental result in optimal control theory, the solutions of this Hamiltonian system can be candidates for the maximal solutions of the optimal control system. We will continue the discussion in this regard to the next sections of the paper.
\subsection{Reduction in optimal control problem on Lie groupoid to its Lie algebroid}
As it is claimed in \cite{Joz, GJ}, every right-invariant control system on a Lie groupoid $G$ reduces to a system on the associated Lie algebroid $AG$.
In other words, equivariant control systems and  optimal control problems on a Lie groupoid $G$ lead naturally to systems and problems on the associated Lie algebroid $AG$.

\begin{example}

 A right-invariant control system on the trivial Lie groupoid 
 $G=M \times  \mathbf{G}\times M$, 
 where $M$ is a smooth manifold and $\mathbf{G}$ is a Lie group, is determined by a vector field on $M$ and a right-invariant vector field on the Lie group $G$.
In other words, a right-invariant control system (\ref{cs}) on $G$  reduces to a control system of the form (\ref{cs2}) on the trivial Lie algebroid 
$$AG= TM + (M\times \mathbf{g}),$$
where $\mathbf{g}$ is the Lie algebra of $\mathbf{G}$.
\end{example}

In general, let we have a right-invariant control system (\ref{cs}) on the Lie groupoid $G.$ We define reduced control system on the Lie algebroid $AG$ of $G$ as follows:

$$f: M\times U \to AG , \quad  f(x, u)= F(1_x, u)= dR_{g^{-1}}F(g, u).$$
Since $F: G\times U \to T^\alpha G$, so for every $g\in G$, $u\in U$  we have  $F(g, u) \in T_g^\alpha G$. On the other hand, $(AG)_x=T_{1_x}^\alpha G.$
so $F(1_x, u)\in T_{1_x}^\alpha G = (AG)_x.$ Therefore $f: M\times U \to AG$, $f(x, u)= F(1_x, u)$ is well-defined.

As a result, on the Lie algebroid $AG$, we have a control system reduced in the control system (\ref{cs}) on the Lie groupoid $G$.

From the above discussion, we conclude that:
\begin{corollary}
Every right-invariant control system and every optimal control problem on a Lie groupoid reduce in its co-adjoint Lie algebroid.
\end{corollary}
So in the latter, we only study control systems and optimal control problems on the co-adjoint Lie algebroids.
\section{Optimal control problem on the co-adjoint Lie algebroid}
Some contents of this section have been taken in \cite{HA1, Mar, HA3}.\\
Let $G\rightrightarrows M$ be a regular Lie groupoid with structural mappings $\alpha, \beta, 1, i, m$. Also let $(AG, M, [., .], \rho, \tau)$ be the Lie algebroid of $G$ over $M$, where $[., .]$ is the bracket on the sections of vector bundle $AG$ over $M$ and $\rho$ is the anchor $\rho: AG \to TM$   and $\tau: AG \to M$  is the bundle map.
 Let us denote the co-adjoint Lie groupoid of the Lie groupoid $G\rightrightarrows M$ by  $\mathcal{G} \rightrightarrows M$ with the structural mappings $\tilde{\alpha}, \tilde{\beta}, \tilde{1}, \tilde{i}$ and $\tilde{m}$. 
Let $(A\mathcal{G}, M, \tilde{[., .]}, \tilde{\rho}, \tilde{\tau})$ be its Lie algebroid.( For more details see \cite{HA1}). As we showed in \cite{HA1}, the co-adjoint orbit $$\mathcal{G}=\{(\mathrm{Ad}_g^{\ast}\xi): g\in G\},$$ for an arbitrary element $\xi$ of $(A^{\ast}I_G)_x$, where $A^{\ast}I_G$ is dual bundle of the isotropy Lie algebroid $A{I_G}$ (see \cite{HA1}).
Also in \cite{HA1} we showed that co-adjoint Lie algebroid associated with the co-adjoint Lie groupoid is 
$$A\mathcal{G}=(T_\xi^{\tilde{\alpha}}\mathcal{G})_{\tilde{1_x}}=(\{\mathrm{ad}_X^{\ast}\xi: X\in T_gG\})_{\mathrm{Ad}_{1_x}^{\ast}}.$$
Moreover, it is proven that for the Lie algebroids $AG$ and $A\mathcal{G}$ we have 
$$\widetilde{[\tilde{X}, \tilde{Y}]}=\mathrm{ad}_{[X, Y]}^\ast \xi,$$
where $\tilde{X}=\mathrm{ad}_X^{\ast}\xi$ and $\tilde{Y}=\mathrm{ad}_Y^{\ast}\xi$  are right-invariant vector field on $\mathcal{G}$ and $ X,  Y$ are right-invariant vector fields on $G$.

\subsection{Poisson structure on dual of the co-adjoint Lie algebroid}
It is well-known fact, there exists  a linear Poisson structure on the dual of any Lie algebroid.(For more details see \cite{DeL}, \cite{Marr}) 

Let  $X$ be a section of $\tau : AG \longrightarrow M, $ the linear function $\hat{X}$  is defined on $A^{\ast} G$ as follows:
$$\hat{X} : A^{\ast} G \longrightarrow \mathbb{R}, $$
 $$\hat{X} (\theta) = \theta ( X ( \tau^{\ast}  (\theta))),$$
where $\theta \in A^{\ast} G$ and  $\tau^{\ast}: A^{\ast} G \longrightarrow M$ is the dual bundle of $\tau : AG\longrightarrow M.$\\

The linear Poisson structure on $A^{\ast} G,$ which is indicated by $\lbrace. , . \rbrace _{A^{\ast} G},$ is characterized by the following conditions:
$$ \lbrace. , . \rbrace _{A^{\ast} G} : C^{\infty} (A^{\ast} G )\times C^{\infty} (A^{\ast} G ) \longrightarrow C^{\infty} (A^{\ast} G ),$$
$$ \lbrace \hat{X} , \hat{Y} \rbrace _{A^{\ast} G} = - [\vert \widehat{ X , Y} \vert], $$ 
$$ \lbrace f \circ \tau^{\ast}  ~,~ \hat{X} \rbrace _{A^{\ast} G} = ( \rho (X) (f) ) \circ \tau^{\ast}, $$
$$\lbrace f \circ \tau^{\ast} ~,~  g \circ \tau^{\ast}  \rbrace _{A^{\ast} G} = 0,$$
where $\tau^{{\ast}} : A^{\ast} G \longrightarrow M,~ f , g \in C^{\infty} (M)$ and $f \circ \tau^{\ast} , g \circ \tau^{\ast} \in C^{\infty} (A^{\ast} G).$\\
Also, linear Poisson bivector on $A^{\ast} G$ is defined by
$$\Pi_{A^{\ast} G} (d \varphi , d \psi ) = \lbrace \varphi , \psi \rbrace _{A^{\ast} G}, $$
where $\varphi , \psi \in C^{\infty} (A^{\ast} G).$\\
Let $H: A^{\ast} G \longrightarrow \mathcal{R}$ be a smooth function on $A^{\ast} G.$ The Hamiltonian vector field $ \mathcal{X}_{H} ^{\Pi_{A^{\ast} G}}$ of $H$ is defined by
\begin{equation}
 \mathcal{X}_{H} ^{\Pi_{A^{\ast} G}} (F) = \lbrace F, H \rbrace _{A^{\ast} G} = \Pi_{A^{\ast} G} ( d F , d H ),
 \label{Hamiltonian}
\end{equation}
where $ F \in C^{\infty} (A^{\ast} G).$
Now, consider $ (A \mathcal{G}, M, \tilde{ \rho}, \tilde{[\vert ~,~ \vert]} )$ as the associated Lie algebroid to the co-adjoint Lie groupoid $\mathcal{G} \rightrightarrows M$. Using what was previously described in the linear Poisson structure on $A^{\ast} G,$  in the following we will show that $A^{\ast} \mathcal{G}$, the dual of $A \mathcal{G},$ has linear Poisson structure.\\
As we know, every section of $\tilde{\tau}: A\mathcal G \to M$ can be written by $\tilde{X}=ad_X^*\xi$, where $X$ is a section of $\tau: AG \to M$. \\
According to above, for a section $\tilde{X}=ad_{X} ^{\ast} \xi$  of  $\Gamma(\tilde{\tau})$ we consider the associated linear function $\hat{\tilde{X}}$ on $A^{\ast} \mathcal{G}$ as follows:
$$\hat{\tilde{X}} : A^{\ast} \mathcal{G} \longrightarrow \mathbb{R},$$

 $$\hat{\tilde{X}} (\delta) = \delta (\tilde{X}(\tilde{\tau}^{{\ast}}(\delta))),$$
 
 where $\delta \in A^{\ast} \mathcal{G}$ and
  $\tilde{\tau}^{\ast} : A^{\ast} \mathcal{G} \longrightarrow M$
   is dual bundle of $\tilde{\tau} : A \mathcal{G} \longrightarrow M.$\\
 In other words, the above formula indicates that
 $$\hat{\tilde{X}} = \widehat{\mathrm{ad}^{\ast} _{X} \xi}.$$
 Also, the linear Poisson structure on $A^{\ast} \mathcal{G}$  can be considered as
 $$ \lbrace. , . \rbrace _{A^{\ast} \mathcal{G}} : C^{\infty} (A^{\ast} \mathcal{G})\times C^{\infty} (A^{\ast} \mathcal{G} ) \longrightarrow C^{\infty} (A^{\ast} \mathcal{G}),$$
$$ \lbrace \hat{\tilde{X}} , \hat{\tilde{Y}} \rbrace _{A^{\ast} \mathcal{G}}  = - [\vert \widehat {\widetilde{ \tilde{X} , \tilde{Y} \vert]}} = - \widehat{\mathrm{ad}_{[\vert X , Y \vert]} ^{\ast} \xi}.$$ 
 It is easy to check that this Poisson structure on $A^{\ast} \mathcal{G}$  satisfies conditions which are mentioned above.\\
For every Hamiltonian $\tilde{H} : A^{\ast} \mathcal{G} \longrightarrow \mathbb{R}$, the Hamiltonian vector field $\mathcal{X}_{\tilde{H}} ^{\Pi_{A^{\ast} \mathcal{G}}}$  on $A^{\ast}\mathcal{G} $ will be considered as equation (\ref{Hamiltonian}).

So on dual of the co-adjoint Lie algebroid we have Poisson structure and its associated Hamiltonian vector field for any function on the dual vector bundle.

\subsection{The prolongation of vector bundle with respect to a Lie algebroid}
Let us consider a general fiber bundle $\mu: E\to M$ over the state space $M$.
The prolongation of $E$ with respect to the Lie algebroid $\left(A, M, \rho, [., .], \tau \right)$ which denoted by  $\mathcal{T}E$ is defined as follows:

Indeed, the prolongation of $E$  with respect to $A$ is the $A$-tangent bundle to $E$. In other words, if we see the Lie algebroid $A$ as a substitute for the tangent to $M$, then the tangent space to $E$ is not a suitable space for describing dynamics on $E,$ because the projection to $M$ of a vector tangent to $E$ is a vector tangent to $M$ and what we want is an element of $A$, the new tangent bundle to $M$.

This $A$-tangent bundle to $E$ is the vector bundle $\mu_1:\mathcal{T} E \to A$ which its each fiber at the point $p\in E_x$ , $x\in M$ is the vector space $$\mathcal{T}_pE=\{(X, V)\in A_x\times T_pE:\rho(X)=T_p\mu(V)\}.$$
Each element of $\mathcal {T}_pE$ we will denote by $(p, X, V)$ where $p\in E$, $X\in A$ and $V\in T_pE.$ 
It is known that $\mathcal{T}E\rightarrow E$ is a Lie algebroid (for more details,  see \cite{Mar} and \cite{DeL}). 

The anchor of this new Lie algebroid is the map $\rho^\prime: \mathcal{T}E \to TE$ given by $\rho^\prime(p, X, V)=V$ and the bracket is defined in terms of projectable sections as follows.\\
A section $Z$ of $\mathcal{T}_pE$ is of form $Z(p)=(\theta(p), V(p))$, where $\theta$ is a section of $A$ along $\mu$ and $V$ is a vector field on $E$. 

\begin{definition}
A section $Z$ is projectable if there exists a section $\sigma$ of $A$ such that $\mu_2 \circ Z= \sigma \circ \mu$, where  $\mu_2: \mathcal{T}E\to A, \quad \mu_2(p, X, V)=X$. It follows that $Z$ is projectable if and only if $\theta=\sigma \circ \mu$, and therefore $Z$ is of the form $Z(p)=\left(\sigma(m), V(p)\right)$, with $m=\mu(p) \in M$.
\end{definition}
The bracket of two projectable sections $Z_1, Z_2$  is defined by 
$$ \left[Z_1, Z_2\right](p)=\left(p, [\sigma_1, \sigma_2](m), [V_1, V_2](p)\right).$$
 For every $p\in E$ it is clear that $[Z_1, Z_2](p)$ is an element of $\mathcal{T}E$.
\begin{definition}
An element of $\mathcal{T}E$ is said to be vertical if it is in the kernel of $\mu_2$, and thus it is of the form $(p, 0, V)$ with $V$ a vertical vector, tangent to $E$ at $p$, i.e. $V$ is in the kernel of $T_p\mu: TE \to TM$.
\end{definition}
Given a local basis $\{e_\alpha\}$ of sections of $A$ and local coordinates $(x^i, \eta^I)$ on $E$ we consider a local base $\{\mathcal{X}_\alpha, \mathcal{V}_l\}$ of sections of $\mathcal{T}E$ as follows 
$$\mathcal{X}_\alpha(p)=\left(p, e_\alpha(x), \rho^i_\alpha (\frac{\partial}{\partial x^i})_p\right), \quad \mathcal{V}_I(p)=\left(p, 0, (\frac{\partial}{\partial \eta^I})_p\right).$$
Thus, any element $Z$ of $\mathcal{T}E$ at $p$ 
$$Z=\left(z^\alpha e_\alpha(m), \rho^i_\alpha z^\alpha\frac{\partial}{\partial x^i}+v^I(\frac{\partial}{\partial \eta^I})_p\right)$$
can be represented as $$ Z=z^\alpha \mathcal{X}_\alpha+v^I\mathcal{V}_I$$
and $(x^i, \eta^I, z^\alpha, v^I)$ are coordinates on $\mathcal{T}E$.
Vertical elements are therefore linear combination of $\{\mathcal{V}_I\}$.\\
The bracket of the elements of this basis are
$$[\mathcal{X}_\alpha, \mathcal{X}_\beta]=C^\gamma_{\alpha\beta} \mathcal{X}_\gamma, \quad[\mathcal{X}_\alpha, \mathcal{V}_J]=0, \quad [\mathcal{V}_I, \mathcal{V}_J]=0$$
and the anchor applied to the section $Z=z^\alpha \mathcal{X}_\alpha+v^I\mathcal{V}_I$ gives the vector field 
$$\rho^\prime(Z)=\rho^i_\alpha z^\alpha\frac{\partial}{\partial x^i}+v^I\frac{\partial}{\partial \eta^I}.$$
The differential of coordinates and elements of the dual basis is given by:
$$dx^i=\rho^i_\alpha\mathcal{X}^\alpha, \quad d\mathcal{X}^\gamma=-\frac{1}{2}C^\gamma_{\alpha\beta} \mathcal{X}^\alpha\wedge \mathcal{X}^\beta.$$
$$d\eta^I=\mathcal{V}^I, \quad d\mathcal{V}^I=0,$$
where $\{\mathcal{X}^\alpha, \mathcal{V}^I\}$ denotes a dual basis and $\rho^i_\alpha$ and $C^\gamma_{\alpha\beta}$ are local functions of the anchor and the bracket of the Lie algebroid $A$ on $M$, respectively. That is 
$$\rho(e_\alpha)=\rho^i_\alpha\frac{\partial}{\partial x^i}, \quad  [e_\alpha, e_\beta]=C^\gamma_{\alpha\beta}e_\gamma.$$
The differential is determined by the relations
$$dx^i=\rho^i_\alpha, \quad de^\alpha=-\frac{1}{2}C^\gamma_{\alpha\beta} e^\alpha\wedge e^\beta.$$
Now as a vector bundle $E$ we consider the dual of the Lie algebroid $A$, i.e. $E=A^\ast$. In this case, on the prolongation of $A^\ast$ that is on $\mathcal{T}A^\ast$ there exists a canonical symplectic structure $\omega_0$. Its definition is analogous to the definition of the canonical symplectic form on the cotangent bundle. First is defined the canonical 1-form $\theta: \mathcal{T}A^\ast \to \mathcal{R}$ by $$\langle\theta, (\eta, X, V)\rangle=\langle\eta, X\rangle.$$ 
Obviously, it vanishes on vertical sections and its coordinates expression is $\theta=\eta^\alpha\mathcal{X}^\alpha.$
The canonical symplectic form is the differential of the canonical 1-form $\omega=-d\theta$, i.e. 
$$\omega=\mathcal{X}^\alpha \wedge \mathcal{V}^\alpha+\frac{1}{2}C_{\alpha\beta}^\gamma\mathcal{X}^\alpha\wedge\mathcal{X}^\beta.$$
As it is known that the structure of Lie algebroid on $A$ is equivalent to the linear Poisson structure on $A^\ast.$
It is easy to see that the Poisson bracket on $A^\ast$  can be expressed in terms of the symplectic form as follows:
Given a function $H\in C^\infty(A^\ast)$ there exists a unique so-called Hamiltonian section $\sigma_H$ of $\mathcal{T}A^\ast$ such that $i_{\sigma_H}\omega=dH$.
Then the Poisson bracket $\{H, G\}$ of two function $H, G$ on $A^\ast$ is given by $$\{H, G\}=-\omega(\sigma_H, \sigma_G).$$
Also, if $H: A^\ast \to \mathcal{R}$ is a function then its differential is $$dH=\rho_\alpha^i\frac{\partial H}{\partial x^i}\mathcal{X}^\alpha+ \frac{\partial H}{\partial\eta^\alpha}\mathcal{V}^\alpha.$$\\
In the case $E=A^\ast\times_M C$, where $C$ is control space $\pi: C\to M$. The coordinates induced to $\mathcal{T}E $ by coordinates $(x^i, \eta_\alpha, u^c)$ on $E=A^\ast\times_M C$ are denoted as $(x^i, \eta_\alpha, u^c, z^\alpha, v^I, v^c)$  and the associated local basis by $(\mathcal{X}_\alpha, \mathcal{V}_I, \mathcal{V}_c)$.\\
The differential is given by:
$$dx^i=\rho^i_\alpha\mathcal{X}^\alpha, \quad  d\mathcal{X}^\gamma=-\frac{1}{2}C^\gamma_{\alpha\beta} \mathcal{X}^\alpha\wedge \mathcal{X}^\beta.$$
$$d\eta^I=\mathcal{V}^I, \quad d\mathcal{V}^I=0.$$
$$du^c=\mathcal{V}^c, \quad d\mathcal{V}^c=0.$$
The differential of any function $H$ on  $E=A^\ast\times_M C$ is 
$$dH=\rho_\alpha^i\frac{\partial H}{\partial x^i}\mathcal{X}^\alpha+ \frac{\partial H}{\partial\eta^\alpha}\mathcal{V}^\alpha+\frac{\partial H}{\partial u^c}\mathcal{V}^c.$$\\
Now let we have a control system $f$ on a Lie algebroid $\tau: A \to M$ with control space $\pi: C\to M$ that is a section of $A$ along $\pi$. 
As we mentioned above, a trajectory of the system $f$ is an integral curve  of the vector field $\rho(f)$  along $\pi$, i.e. a trajectory is a solution of the following equation
$$\dot{x(t)}=\rho\left(f(c(t))\right),$$  where $x(t)=\pi(c(t)).$\\
An optimal problem on the Lie algebroid $A$ is minimizing  integral of a cost function $L\in C^\infty(C)$ over the set of trajectories of the system which satisfies some boundary conditions. To solve this problem we usually consider a Hamiltonian function $H\in C^\infty (A^\ast\times_M C)$ by $H(\eta, c)=\langle\eta, f(c)\rangle- L(c)$ and the associated Hamiltonian control system $f_H$- a section, defined on a subset of $\mathcal{T}A^\ast$ along $pr_1:A^\ast\times_M C\to A^\ast$ by means of the symplectic equation $$i_{f_H}\omega=dH.$$
The critical trajectories are the integral curves of the vector field $\rho^\prime(f_H)$. The solutions of the optimal control problem are among the critical trajectories. (see \cite{Mar} for details)

Now, consider the  Lie algebroid $ ( \tau : AG \longrightarrow M, \rho, [\vert~,~\vert] ) $ associated to Lie groupoid $G \rightrightarrows M$ with non-zero bracket. 
Furthermore, suppose that $(x^{i})$ are local coordinates on $M$ and $\lbrace e_{\alpha} \rbrace$ is a local basis of sections for $AG.$ 
Consider local functions $\rho_{\alpha} ^{i}, C_{\alpha \beta} ^{\gamma}$ on $M$ which are called structure functions of Lie algebroid $AG$. Also, consider co-adjoint Lie algebroid $( \tilde{\tau}: A\mathcal{G} \longrightarrow M, \tilde{\rho}, \widetilde{[\vert ~,~ \vert]} )$  of the co-adjoint Lie groupoid $\mathcal{G} := O(\xi) \rightrightarrows M.$  

As it is shown in \cite{HA2}, that the basis of sections for co-adjoint Lie algebroid $A\mathcal{G}$, for every $\xi\neq 0$ are as follows:
$$\tilde{e}_{\alpha} = \mathrm{ad}^{\ast} _{e_{\alpha}} \xi.$$
Also, it is proven in \cite{HA2} that  the structure functions of the Lie algebroids $AG$ and $A \mathcal{G}$ are equal, i.e. if  $\rho_{\alpha} ^{i}, C_{\alpha \beta} ^{\gamma}$ are structure functions of the Lie algebroid $AG$ and 
$\tilde{\rho}_{\alpha}^{i}, \tilde{C}_{\alpha \beta}^{\gamma}$ are structure functions of Lie algebroid $A \mathcal{G},$ then 
$$\tilde{\rho}_{\alpha}^{i}=\rho_{\alpha} ^{i}, \quad \tilde{C}_{\alpha \beta}^{\gamma} =C_{\alpha \beta} ^{\gamma}. $$

Now for local coordinates $(x^i, y^\alpha)$ on $A\mathcal{G}$ associated to the base$ \{\tilde{e}_\alpha\}$ of sections of $A\mathcal{G}$ and coordinates $(x^i, \eta^\alpha)$ on $A^\ast\mathcal{G}$ we define a local base $\{\tilde{\mathcal{X}}_\alpha , \tilde{\mathcal{V}}_\alpha\}$ of sections for the  prolongation of $\mathcal{T}A^\ast\mathcal{G}$ with respect to the co-adjoint Lie algebroid $A\mathcal{G}$ by   

$$\tilde{\mathcal{X}}_\alpha(p)=\left(p, \mathrm{ad}^\ast_{e_\alpha}\xi(x), \rho^i_\alpha (\frac{\partial}{\partial x^i})_p\right), \quad \tilde{\mathcal{V}}_\alpha(p)=\left(p, 0, (\frac{\partial}{\partial \eta^\alpha})_p\right),$$
where $p\in (A^\ast\mathcal{G})_x$.\\
The canonical symplectic form in these coordinates is as follows:
$$\omega=\tilde{\mathcal{X}}^\alpha \wedge\tilde{\mathcal{V}}^\alpha+\frac{1}{2}C_{\alpha\beta}^\gamma\tilde{\mathcal{X}}^\alpha \wedge\tilde{\mathcal{X}}^\beta.$$
As $\{\tilde{\mathcal{X}}_\alpha, \tilde{\mathcal{V}}_\alpha\}$ is a basis for sections of $A\mathcal{G}$ so the local expression of the section $f_H$ is 
 $$f_H= \lambda^\alpha\tilde{\mathcal{X}}_\alpha+ \mu^\alpha\tilde{\mathcal{V}}_\alpha.$$
Then we obtain
$$i_{f_H}\omega=\lambda^\alpha \tilde{\mathcal{V}}^\alpha-\left(\mu^\alpha+\eta_\gamma C^\gamma_{\alpha\beta}\lambda^\beta\right) \tilde{\mathcal{X}}^\alpha.$$
By definition $dH$ and substitute it in above relations, we have
$$f_H= \frac{\partial H}{\partial \eta_\alpha}\tilde{\mathcal{X}}_\alpha 
\left(\rho^i_\alpha\frac{\partial H}{\partial x^i}+
 \eta_{\gamma}C^{\gamma}_{\alpha \beta}\frac{\partial H}{\partial \eta_{\beta}}\right)
 \tilde{\mathcal{V}}^\alpha$$ 
defined on the subset $$\frac{\partial H}{\partial u^c}=0$$
and the coordinate expression of the vector field $\rho^\prime(f_H)$  is
$$\rho^\prime(f_H)=  \rho^i_\alpha \frac{\partial H}{\partial \eta_\alpha}
\frac{\partial}{\partial x^i}- \left(\rho^i_\alpha\frac{\partial H}{\partial x^i}+
 \eta_\gamma C^\gamma_{\alpha\beta}\frac{\partial H}{\partial \eta_\beta}\right)\frac{\partial}{\partial \eta_\alpha}$$ 
so as a main result of this work,  we obtain that the critical trajectories are the solution of the following differential equations:\\
\begin{eqnarray}{}
\label{critical tra}
\dot{x^i}&=&\rho^i_\alpha \frac{\partial H}{\partial \eta_\alpha},\nonumber\\
\dot{\eta_\alpha}&=&- \left(\rho^i_\alpha\frac{\partial H}{\partial x^i}+ \eta_\gamma C^\gamma_{\alpha\beta}\frac{\partial H}{\partial \eta_\beta}\right), \\
0&=&\frac{\partial H}{\partial u^c}.\nonumber
\end{eqnarray}

\section{Example}

Let $\mathbf{G}$ be a Lie group and $M$ be a manifold. Consider the trivial Lie groupoid  $G := M \times \mathbf{G} \times M \rightrightarrows M.$ As described in \cite{Mac1}, the Lie algebroid associated to the trivial Lie algebroid is $A G = TM \oplus ( M \times \mathbf{g} ).$ The anchor $\rho :  TM \oplus ( M \times \mathbf{g} ) \longrightarrow TM$ is the projection $X \oplus V \longmapsto X$ and the Lie bracket on the sections of $ TM \oplus ( M \times \mathbf{g} )$ is given by
$[\vert X \oplus V, Y \oplus W \vert] = [X,Y] \oplus \lbrace X(W) - Y (V) + [V,W]\rbrace.$\\

The co-adjoint Lie groupoid  associated to the trivial Lie groupoid is $\mathcal{G}_\xi:= M \times O(\xi^{\prime}) \rightrightarrows M,$ where $O(\xi^{\prime}) $ is the orbit of co-adjoint action of Lie group and 
$\xi^{\prime} \in g^\ast$  for which $G_{\xi^{\prime}}=\{a\in \mathbf{G}:Ad_a^\ast \xi^{\prime}=\xi^\prime\}$ is the normal Lie subgroup of $\mathbf{G}$. The Lie algebroid of the co-adjoint Lie groupoid is $A \mathcal{G}:= M \times T_{\xi^{\prime}} O(\xi^{\prime}).$ The anchor $\tilde{\rho}$ is given by
\begin{eqnarray}
\label{anchor}
&~&\tilde{\rho} : M \times T_{\xi^{\prime}} O (\xi^{\prime} )  \longrightarrow  TM, \nonumber  \\
&~&\tilde{\rho}  (x , \mathrm{ad}_{V} ^{\ast}  \xi^{\prime} ) (p) = X (p),
\end{eqnarray}
where $X \in \Gamma (TM)= \mathbf{X}(M)$ is equal to $\dot{p}(0)$, $p(t)=\beta(\gamma(t))\in M$, $\gamma(t)=(p(t), \mathrm{Ad}_{a} ^{\ast}  \xi^{\prime}, p(t) )\in \mathcal{G}$, $(\frac{d}{dt})_{t=0}(\gamma(t))=(x, \mathrm{ad}_{V} ^{\ast}  \xi^{\prime})\in A \mathcal{G}$ and $ p(0)=p \in M.$\\ 
The Lie bracket on the space of sections of $ A \mathcal{G} $ is
\begin{eqnarray}
\label{bracket}
[\vert  \mathrm{ad}_{V} ^{\ast}  \xi^{\prime} , \mathrm{ad}_{W} ^{\ast}  \xi^{\prime}  \vert] = \mathrm{ad}^{\ast} _{[V,W]} \xi^{\prime}, \nonumber
\end{eqnarray}
for every $V^\prime = \mathrm{ad}_{V} ^{\ast}  \xi^{\prime} ,~ W^\prime = \mathrm{ad}_{W} ^{\ast}  \xi^{\prime}  \in \Gamma ( M \times T_{\xi^{\prime}} O (\xi^{\prime} ) )$, where $V, W \in \mathbf{g}$ .\\
As we know, if we assume that $ \hat{V^\prime} , \hat{W^\prime} \in C^{\infty} ( T^{\ast} _{\xi^{\prime}} O(\xi^{\prime} )), $ then $T^{\ast} _{\xi^{\prime}} O(\xi^{\prime} )$ carries the Kirillov-Kostant bracket as follows:
$$\lbrace \hat{V^\prime} , \hat{W^\prime} \rbrace ( \lambda ) = \langle \lambda , [\vert V^\prime , W^\prime \vert ] \rangle. $$
Now, consider vector bundle $\tau : M \times T_{\xi^{\prime}} O(\xi^{\prime} ) \longrightarrow M,$ which is projection over the first factor, and its dual $\tau^{\ast} : M \times T^{\ast} _{\xi^{\prime}} O(\xi^{\prime} ) \longrightarrow M.$ Let $\Sigma^{\prime} :=(p,V^\prime) = (p, \mathrm{ad}_{V} ^{\ast} \xi^{\prime}) \in M \times T _{\xi^{\prime}} O(\xi^{\prime} )$ and $\delta = (p, \lambda ) \in M \times T^{\ast} _{\xi^{\prime}} O(\xi^{\prime} ), $ so the linear function $\hat{\Sigma^{\prime}}$ on $M \times T^{\ast} _{\xi^{\prime}} O(\xi^{\prime} )$ will be as follows:
$$\hat{\Sigma^{\prime}} : M \times T^{\ast} _{\xi^{\prime}} O(\xi^{\prime} ) \longrightarrow \mathbb{R},$$
 $$\hat{\Sigma^{\prime}} (\delta) = \delta ( \Sigma^{\prime} ( \tau^{\ast} (\delta))= \langle \lambda , V^{\prime} \rangle.$$
In other word we have $\hat{\Sigma^{\prime}} =(p, \hat{V^{\prime}} ).$\\
Now, we try to clear relation between $C^{\infty} ( M \times T^{\ast} _{\xi^{\prime}} O(\xi^{\prime} ))$ and $C^{\infty} ( T^{\ast} _{\xi^{\prime}} O(\xi^{\prime} )).$ Let $ V^{\prime} \in C^{\infty} ( T^{\ast} _{\xi^{\prime}} O(\xi^{\prime} ))$ be a vector field on $O(\xi^{\prime} ).$ We define
\begin{eqnarray}
\hat{V^{\prime}}:  T^{\ast} _{\xi^{\prime}} O(\xi^{\prime} ) ) \longrightarrow \mathbb{R}, \nonumber   \\
\hat{V^{\prime}} (\lambda ) = \langle \lambda , V^{\prime} \rangle = \lambda (V^{\prime}).\nonumber 
\end{eqnarray}
$\hat{V^{\prime}}$ is linear function on $ T^{\ast} _{\xi^{\prime}} O(\xi^{\prime} ).$\\
Now, Kirillov-Kostant bracket on $C^{\infty} ( T^{\ast} _{\xi^{\prime}} O(\xi^{\prime} )) $ is as follows:
 $$ \lbrace. , . \rbrace _{ K.K } : C^{\infty} ( T^{\ast} _{\xi^{\prime}} O(\xi^{\prime} ))\times C^{\infty} ( T^{\ast} _{\xi^{\prime}} O(\xi^{\prime} ) ) \longrightarrow C^{\infty} (T^{\ast} _{\xi^{\prime}} O(\xi^{\prime} )),$$
$$( \hat{V^{\prime}} , \hat{W^{\prime}} ) \longmapsto \lbrace \hat{V^{\prime}} , \hat{W^{\prime}} \rbrace_{K.K},$$
$$\lbrace \hat{V^{\prime}} , \hat{W^{\prime}} \rbrace_{K.K}(\lambda) = - \langle \lambda , [\vert V^{\prime} , W^{\prime} \vert ] \rangle.$$
As we know $ T _{\xi^{\prime}} O(\xi^{\prime} ) = \lbrace \mathrm{ad}_{V} ^{\ast}  \xi^{\prime} ~\vert~ V \in \mathbf{g} \rbrace,$ so let $V^{\prime} = \mathrm{ad}_{V} ^{\ast}  \xi^{\prime} \in  T _{\xi^{\prime}} O(\xi^{\prime} ).$ Also by well known fact for finite dimensional vector space
$$ T^{\ast \ast} _{\xi^{\prime}} O(\xi^{\prime} ) \cong T _{\xi^{\prime}} O(\xi^{\prime} ),$$
 one can consider $V^{\prime} \in T^{\ast \ast} _{\xi^{\prime}} O(\xi^{\prime} ), $ i.e. $V^{\prime} :  T^{\ast} _{\xi^{\prime}} O(\xi^{\prime} ) \longrightarrow \mathcal{R}$ is linear functional, in other words, $V^{\prime} \in C^{\infty} ( T^{\ast} _{\xi^{\prime}} O(\xi^{\prime} ) ).$ So $T _{\xi^{\prime}} O(\xi^{\prime} ) \subset C^{\infty} (T^{\ast} _{\xi^{\prime}} O(\xi^{\prime} ) ),$ therefore we can take $\hat{V^{\prime}} = V^{\prime}$ and for every $\lambda \in T^{\ast} _{\xi^{\prime}} O(\xi^{\prime} )$ we have
$$\hat{V^{\prime}} (\lambda ) = V^{\prime} (\lambda ) = \langle \lambda , V^{\prime} \rangle$$
or equivalently 
$$\langle \lambda , V^{\prime} \rangle =( \mathrm{ad}_{V} ^{\ast} \xi^{\prime}) (\lambda ).$$
Now, we rewrite the bracket 
 $$ \left\lbrace. , . \right\rbrace _{ K.K } : C^{\infty} ( T^{\ast} _{\xi^{\prime}} O(\xi^{\prime} ))\times C^{\infty} ( T^{\ast} _{\xi^{\prime}} O(\xi^{\prime} ) ) \longrightarrow C^{\infty} (T^{\ast} _{\xi^{\prime}} O(\xi^{\prime} ))$$
 as follows:
 \begin{eqnarray}
  \left\lbrace \hat{V^{\prime}} , \hat{W^{\prime}} \right\rbrace_{K.K} (\lambda ) &=& - \langle \lambda , [\vert V^{\prime} , W^{\prime} \vert] \rangle, \nonumber  \\
  &=& - [\vert V^{\prime} , W^{\prime} \vert ] (\lambda),  \nonumber \\
  &=& - [ \vert \mathrm{ad}_{V} ^{\ast} \xi^{\prime},  \mathrm{ad}_{W} ^{\ast} \xi^{\prime} \vert ] (\lambda ),  \nonumber  \\
  &=& - \mathrm{ad}_{[V , W]} ^{\ast} \xi^{\prime} (\lambda ), \nonumber
 \end{eqnarray}
 as a result we obtain
 $$\lbrace \hat{V^{\prime}} , \hat{W^{\prime}} \rbrace_{K.K} = - \widehat{\mathrm{ad}_{[V , W]} ^{\ast} \xi^{\prime}}=-\mathrm{ad}_{[V , W]} ^{\ast} \xi^{\prime}.$$
So, according to the subsection 3.3 in \cite{HA1}, one can easily check that the first property of the linear Poisson structure of functions
 $$ \lbrace. , . \rbrace _{ A^{\ast} \mathcal{G} } : C^{\infty} (M \times T^{\ast} _{\xi^{\prime}} O(\xi^{\prime} ))\times C^{\infty} (M \times T^{\ast} _{\xi^{\prime}} O(\xi^{\prime} ) ) \longrightarrow C^{\infty} (M \times T^{\ast} _{\xi^{\prime}} O(\xi^{\prime} )) $$
 will be as follows:
\begin{eqnarray}
\label{kk}
 \lbrace \hat{\Sigma^{\prime} _{1}} , \hat{\Sigma^{\prime} _{2}} \rbrace _{A^{\ast} \mathcal{G}} (\delta ) = \lbrace \hat{V^{\prime}} , \hat{W^{\prime}} \rbrace_{K.K} (\lambda ).\nonumber
\end{eqnarray}
where $ \hat{\Sigma^{\prime} _{1}} = (p, \hat{V^{\prime}}) , \hat{\Sigma^{\prime} _{2}} = (p, \hat{W^{\prime}})$ and $\delta = (p, \lambda ) \in M \times T^{\ast} _{\xi^{\prime}} O(\xi^{\prime} ).$ So, we have a well known Poisson structure on $M \times T^{\ast} _{\xi^{\prime}} O(\xi^{\prime} ).$\\
 Also, according to the subsection 3.3 in \cite{HA1} and  equation (\ref{anchor}), the second property of the linear Poisson structure on $M \times T^{\ast} _{\xi^{\prime}} O(\xi^{\prime} ) $ is 
 $$ \lbrace f \circ \tau^{\ast}  , \hat{\Sigma^{\prime}}  \rbrace_{A^{\ast} \mathcal{G}} (\delta)= \left(\tilde{\rho} (\Sigma^{\prime} ) (f) \right)\circ \tau^{\ast} (\delta) = X (f(p)),$$  
 where $\Sigma^{\prime} = \mathrm{ad}_{\Sigma} ^{\ast} \xi^{\prime} \in \Gamma (A \mathcal{G} )$ and $\Sigma = X \oplus V \in AG. $

 Moreover, as we mentioned in the subsection 3.3 in \cite{HA1}, the third feature of linear Poisson structure on $M \times T^{\ast} _{\xi^{\prime}} O(\xi^{\prime} )$ easily deduced based on the linear Poisson structure on $A^{\ast}G,$ i.e.
\begin{eqnarray}
\lbrace f \circ \tau^{\ast} , g \circ \tau^{\ast} \rbrace_{A^{\ast} \mathcal{G}} = 0.  \nonumber  
\end{eqnarray}
 
Now, suppose that $H: M \times T^{\ast} _{\xi^{\prime}} O(\xi^{\prime} ) \longrightarrow \mathcal{R}$ be function which we define it as $ H =(p, h)$ where $h: T^{\ast} _{\xi^{\prime}} O(\xi^{\prime} ) \longrightarrow \mathcal{R}$ is a Hamiltonian function.
In the following, we will show that $H$ is a Hamiltonian function on $M \times T^{\ast} _{\xi^{\prime}} O(\xi^{\prime} ).$ In order to reach this result, we need to express some fundamental information which are related to Hamiltonian mechanics on cotangent bundles and Lie algebroids.\\

Consider $O(\xi^{\prime} )$ as a smooth manifold and let $T^{\ast} O(\xi^{\prime} )$ be its cotangent bundle. Suppose that $\delta = (p,\lambda ) \in T^{\ast} O(\xi^{\prime} )$ and $X_{\delta} \in T_{\delta} ( T^{\ast} O(\xi^{\prime} ) ).$ As we know, the  Liouville form on $  T^{\ast} O(\xi^{\prime} ) $ is the 1-form  $\theta$ such that
\begin{eqnarray}
\label{big teta}
\theta (X_{\delta} ) = \lambda \left(T \pi_{O(\xi^{\prime} )} (X_{\delta} )\right).\nonumber
\end{eqnarray}
where $\pi_{O(\xi^{\prime} )} : T^{\ast} O(\xi^{\prime} ) \longrightarrow O(\xi^{\prime}) ;  ~ (p, \lambda )\longmapsto p $ is canonical projection. Moreover, the 2-form 
\begin{eqnarray}
\omega = d \theta 
\label{omega}
\end{eqnarray}
is canonical symplectic form on $T^{\ast} O(\xi^{\prime} ),$ \\

Furthermore, a vector field $\mathcal{X}, $ where $\mathcal{X} \in \Gamma ( T ( T^{\ast} O(\xi^{\prime} ) ) ), $ is called Hamiltonian vector field if there is a $h \in C^{\infty} (T^{\ast}_{\xi^\prime} O(\xi^{\prime} ) )$ such that $i_\mathcal{X} \omega = dh.$ Let $h : T^{\ast}_{\xi^\prime} O(\xi^{\prime} ) \longrightarrow \mathcal{R}$ be a Hamiltonian function and $\mathcal{X}_{h}$ be Hamiltonian vector field associated to Hamiltonian function $h.$ Moreover,
\begin{eqnarray}
\label{kk}
\mathcal{X}_{h} (f ) = \Pi (d f , d h ) = \lbrace f, h \rbrace_{K.K}, \nonumber
\end{eqnarray}
where $\Pi$ is Poisson 2-vector on $ T^{\ast} O(\xi^{\prime} ) $,  $f \in C^{\infty} (T^{\ast}  O(\xi^{\prime} ))$ and $\lbrace f, h \rbrace_{K.K}$ is the Kirillov-Kostant bracket on   $ C^{\infty} (T^{\ast}_{\xi^\prime} O(\xi^{\prime} ) )$(see \cite{HA1} for more details).\\

In local coordinates $\{\eta_\alpha\}$ for $ T^{\ast} O(\xi^{\prime} ) $, the Poisson 2-vector $\Pi$ is $$\Pi=-\frac{1}{2} C^\gamma_{\alpha\beta}\eta_\gamma\frac{\partial}{\partial \eta_\alpha} \wedge \frac{\partial}{\partial \eta_\beta}$$ and 
Hamiltonian vector field $\mathcal{X}_{h}$  associated to Hamiltonian function $h$  is $$\mathcal{X}_{h}=-C^\gamma_{\alpha\beta}\eta_\gamma \frac{\partial h}{\partial \eta_\beta}\frac{\partial}{\partial \eta_\alpha}.$$ 
So Hamiltonian equations are
\begin{equation}
\label{HEh}
\dot{\eta_\alpha}=-C^\gamma_{\alpha\beta}\eta_\gamma \frac{\partial h}{\partial \eta_\beta}.
\end{equation}

Now, consider the prolongation $\mathcal{T}A^{\ast} \mathcal{G}$
\begin{eqnarray}
\mathcal{T}A^{\ast} \mathcal{G} &=& \big\lbrace  (\delta , a , \vartheta_{\delta} ) \in M \times T^{\ast}  O(\xi^{\prime} ) \times T O(\xi^{\prime} ) \times T_{\delta} (T^{\ast}  O(\xi^{\prime} )) ~ \big\vert ,  \nonumber   \\ 
 &~& \tilde{\rho }( a ) = T \tau^{\ast} (\vartheta_{\delta}),~ \vartheta_{\delta} \in T_{\delta} A^{\ast} \mathcal{G},~ \tau^{\ast} (\delta) = \tau (a) \big\rbrace,  \nonumber
\end{eqnarray}
where $\delta =(p, \lambda) \in M \times T^{\ast} O(\xi^{\prime} )$ and $ a = (p , V^{\prime} )\in M \times T O(\xi^{\prime} ).$ 

The vector bundle $\tau^\prime : \mathcal{T} A^{\ast} \mathcal{G}\longrightarrow A^{\ast} \mathcal{G}$ has Lie algebroid structure $( \rho^\prime , [\vert ~,~ \vert]^\prime) $ such that
\begin{enumerate}
\item The anchor $\rho^\prime : \mathcal{T} A^{\ast} \mathcal{G} \longrightarrow T A^{\ast} \mathcal{G}$ is projection onto the third factor, 
\begin{eqnarray}
\rho^\prime (\delta, a, \vartheta_{\delta} ) = \vartheta_{\delta}. \nonumber 
\end{eqnarray}
\item  A section $\tilde{\Sigma} \in \Gamma (\tau^\prime)$ is projectable if there exists a section $\Sigma$ of $\tau : A \mathcal{G} \longrightarrow M $ and a vector field $ \mathcal{X} \in \mathbf{X} (A^{\ast} \mathcal{G} )$ which is $\tau$-projectable to the vector field $\tilde{\rho} (\Sigma )$ on $M,$ such that $\tilde{\Sigma}( \delta ) = ( \delta, \Sigma (p), \mathcal{X}(\delta) )$ for all $\delta\in A^{\ast} \mathcal{G}.$ We use the notation $\tilde{\Sigma} \equiv (\Sigma , \mathcal{X} ).$ Then the bracket of two projectable sections $\tilde{\Sigma_{1}}$ and $\tilde{\Sigma_{2}}$ is given by
$$[\vert \tilde{\Sigma}_{1}, \tilde{\Sigma}_{2} \vert] (\delta) 
= \left( \delta, [\vert \Sigma _{1}, \Sigma _{2} \vert] (\tau (\delta)) , [ \mathcal{X}_{1} , \mathcal{X}_{2} ] (\delta)) \right).$$
\end{enumerate}
 In \cite{HA1} it is shown that for every  Liouville 1-form $\theta$  on $T^{\ast} O(\xi^{\prime} )$ there exists a Liouville section $\Theta \in \Gamma ( (\mathcal{T}A^{\ast} \mathcal{G} )^{\ast} )$ such that : 
 \begin{eqnarray}
 \label{teta}
 \Theta = (p , \theta ).
 \end{eqnarray}

Furthermore, according to equations (\ref{omega}) and (\ref{teta}), the canonical symplectic section $ \Omega $  will be defined as follows:
\begin{eqnarray}
\Omega = - d \Theta = (p, \omega), \nonumber
\label{Omg}
\end{eqnarray}
where $\omega$ is canonical symplectic 2-form on $T^{\ast} O(\xi^{\prime} ).$\\

Let $H : M \times T^{\ast} _{\xi^{\prime}} O(\xi^{\prime} ) \longrightarrow \mathcal{R}$ be a Hamiltonian function, $\Omega$ be symplectic section and $d H \in \Gamma ( ( \mathcal{T} A^{\ast} \mathcal{G} )^{\ast}).$ Then, by definition, there exists the unique Hamiltonian section $\mu_{H} \in \Gamma ( \mathcal{T} A^{\ast} \mathcal{G} ) $ satisfying
 $$i_{\mu_{H}} \Omega = d H.$$
 
In the following lemma, we will show the correspondence between Hamiltonian sections associated to Lie algebroid $ M \times T_{\xi^{\prime}} O(\xi^{\prime} )$ and tangent space $ T_{\xi^{\prime}} O(\xi^{\prime} ).$
 \begin{lemma}
 Consider Hamiltonian function $ H= (p , h): M \times T^{\ast} _{\xi^{\prime}} O(\xi^{\prime} ) \longrightarrow \mathbb{R} $ where $h: T^{\ast} _{\xi^{\prime}} O(\xi^{\prime} ) \longrightarrow \mathbb{R} $ is Hamiltonian function defined on $T^{\ast} _{\xi^{\prime}} O(\xi^{\prime} ).$ Let $\mathcal{X}_{h}$ be Hamiltonian vector field of $h.$ Then the Hamiltonian section of H will be as follows:
 $$\mu_{H} = \left(p , \mathcal{X}_{h} \right).$$
 \end{lemma}
 
{ \bf Proof:}
 
 Let $Y^{\prime} = (p , Y) $ be a vector field on $M \times T^{\ast} _{\xi^{\prime}} O(\xi^{\prime} ),$ Then
 \begin{eqnarray}
 d H (Y^{\prime} ) &=& (p , dh) (p, Y) \nonumber  \\
 &=&  i_{\mathcal{X}_{h}} \omega (Y)  \nonumber  \\
 &=&  \omega (\mathcal{X}_{h} , Y)   \nonumber  \\
 &=&  (p , \omega ) ( (p, \mathcal{X}_{h} ) (p, Y) )\nonumber  \\
 &=&  \Omega ( (p, \mathcal{X}_{h}) , Y^{\prime} ) \nonumber  \\
 &=& i_{(p, \mathcal{X}_{h})} \Omega (Y^{\prime}).\nonumber
 \end{eqnarray}
 So, if $\mathcal{X}_{h}$ is the Hamiltonian vector field associated to Hamiltonian function $h,$ then, according to what was presented above, and since the Hamiltonian vector field associated to Hamiltonian function $h$ and Hamiltonian section associated to Hamiltonian function are unique, we conclude that the $\mu_{H} = (p, \mathcal{X}_{h} )$ is Hamiltonian section associated to Hamiltonian function $H,$ and vice versa.

 Furthermore, $\tilde{\rho} ( \mu_{H} )$ is Hamiltonian vector field of $H$ with respect to the linear Poisson structure $\Pi_{M \times T^{\ast} _{\xi^{\prime}} O(\xi^{\prime} )}$ on $M \times T^{\ast} _{\xi^{\prime}} O(\xi^{\prime} ).$
 So, according to definition of the anchor $\tilde{\rho}$, we have that
 $$\tilde{\rho} ( \mu_{H} ) = \mathcal{X}_{H} \in \mathbf{X} ( M \times T^{\ast} _{\xi^{\prime}} O(\xi^{\prime} )).$$
  We denote by $\mathcal{X}_{H} ^{\Pi_{ M \times T^{\ast} _{\xi^{\prime}} O(\xi^{\prime} )}} $ the Hamiltonian vector field of $H$ with respect to the linear Poisson structure $\Pi_{ M \times T^{\ast} _{\xi^{\prime}} O(\xi^{\prime} )}$ on $M \times T^{\ast} _{\xi^{\prime}} O(\xi^{\prime} ).$ \\
Note that by using the equations (\ref{Hamiltonian}) and (\ref{kk}), we have actually proved that:
 $$ \mathcal{X}_{H} ^{\Pi_{M \times T^{\ast} _{\xi^{\prime}} O(\xi^{\prime} )}} (F) =\mathcal{X}_{h} (f ),$$
where $ F = (p , f ) \in C^{\infty} (M \times T^{\ast} _{\xi^{\prime}} O(\xi^{\prime} )) $ and $ f \in C^{\infty} (T^{\ast} _{\xi^{\prime}} O(\xi^{\prime} ))$(see \cite{HA1} for more details).\\

Let $(x^{i})$ be local coordinates on open subset $U$ of $M,~ \lbrace e_{\alpha} \rbrace$ is local basis of sections for $A,$ we have that 

\begin{equation}
\Pi_{A^{\ast}} = \rho_{\alpha} ^{i}  \frac{\partial}{\partial x^{i}} \wedge \frac{\partial}{\partial y_{\alpha}} - \frac{1}{2} C_{\alpha \beta} ^{\gamma} y_{\gamma} \frac{\partial}{\partial y_{\alpha}} \wedge \frac{\partial}{\partial y_{\beta}}
\label{Poisson}
\end{equation}
where $(x^{i} , y_{\alpha})$ are the corresponding local coordinates on $A^{\ast}$ and $\rho_{\alpha} ^{i}, C_{\alpha \beta} ^{\gamma}$ are the local structure functions of $A$ with respect to the coordinates $(x^{i} )$ and basis $  \lbrace e_{\alpha} \rbrace .$\\

For Hamiltonian function $H: A^{\ast} \longrightarrow \mathcal{R},$  the Hamiltonian vector field associated to $\Pi_{A^{\ast}}$ is as follows:
\begin{equation}
 \mathcal{X}_{H} ^{\Pi_{A^{\ast} }} (F) = \Pi_{A^{\ast} } ( d F , d H )=\{F, H\}_{A^\ast\mathcal{G}},
 \label{HVF}
\end{equation}
where $ F \in C^{\infty} (A^{\ast}).$
From equations (\ref{Poisson}) and (\ref{HVF}) it follows that the local expression of $ \mathcal{X}_{H} ^{\Pi_{A^{\ast} }}$ is: 
\begin{equation}
 \mathcal{X}_{H} ^{\Pi_{A^{\ast} }} = \frac{\partial H}{\partial y_{\alpha}}  \rho_{\alpha} ^{i}  \frac{\partial}{\partial x^{i}}  - \left( \frac{\partial H}{\partial x^{i}} \rho_{\alpha} ^{i} +  \frac{\partial H}{\partial y_{\beta}} C_{\alpha \beta} ^{\gamma} y_{\gamma} \right) \frac{\partial}{\partial y_{\alpha}}.
 \label{LHVF}
\end{equation}
So, the Hamiltonian equations are 
\begin{equation}
\label{HEA}
 \frac{d x^{i}}{d t}= \frac{\partial H}{\partial y_{\alpha}}  \rho_{\alpha} ^{i}, \quad\frac{d y_{\alpha}}{d t} = - \left( \frac{\partial H}{\partial x^{i}} \rho_{\alpha} ^{i} +  \frac{\partial H}{\partial y_{\beta}} C_{\alpha \beta} ^{\gamma} y_{\gamma} \right).
\end{equation}

Now we come back to our main example, the co-adjoint Lie groupoid $\mathcal{G}\rightrightarrows M$ of the trivial Lie groupoid $G=M\times \mathbf{G}\times M \rightrightarrows M$.
We discussed in full detail in \cite{HA1} that for co-adjoint Lie algebroid $ (A \mathcal{G}, \tilde{\rho}, [\vert ~, ~ \vert]^{\prime} ),$ its dual bundle, $A^{\ast} \mathcal{G}$ has linear Poisson structure.

Let $(\eta_{\alpha})$ be local coordinates on $T^{\ast} O(\xi^{\prime}) $ and $(x^{i})$ be local coordinates on $M.$ So using the Hamiltonian equations (\ref{HEA}), and similar to what was stated in equation (\ref{LHVF}), and according to that $\rho^i_\alpha=1$ for the Lie algebroid $AG$ and as well for $A\mathcal{G}$, the Hamiltonian vector field $\mathcal{X}_{H} \in \mathbf{X} ( M \times T^{\ast} _{\xi^{\prime}} O(\xi^{\prime} ))$ for Hamiltonian $H: M \times T_{\xi^{\prime}} ^{\ast} O(\xi^{\prime})  \longrightarrow \mathcal{R}$ is as follows:
\begin{equation}
\mathcal{X}_{H} ^{\Pi_{A^{\ast}} \mathcal{G}} =
 \frac{\partial H}{\partial \eta^{\alpha}} \frac{\partial }{\partial x^{i}} 
  - \left(\frac{\partial H}{\partial x^i} + 
   C_{\alpha \beta} ^{\gamma} \eta_{\gamma} \frac{\partial H}{\partial \eta_{\alpha}}\right) \frac{\partial} {\partial  \eta_{\alpha}}.\nonumber
\end{equation}
  Thus, the corresponding Hamiltonian equations are as follows:
 \begin{equation}
 \label{HE1}
 \frac{d x^{i}}{dt} = \frac{\partial H}{\partial \eta^{\alpha}}, \quad   \frac{d \eta_{\alpha}}{dt} =
  - \left(\frac{\partial H}{\partial x^i} + C_{\alpha \beta} ^{\gamma} \eta_{\gamma}\frac{\partial H}{\partial \eta_{\alpha}}\right).
 \end{equation}
 Moreover, if $(x^{i}) $ are local coordinates on $M, \lbrace e^{\prime} _{\alpha} \rbrace $ is the local basis of $\Gamma (A \mathcal{G} )$ and $(x^{i}, \eta_{\alpha})$ are corresponding coordinates on $ A^{\ast} \mathcal{G} = M \times T^{\ast} _{\xi^{\prime}} O(\xi^{\prime} )$ then the local expression of $\Pi_{A^{\ast} \mathcal{G}} $ will be as follows:
 \begin{equation}
 \label{C P}
  \Pi_{A^{\ast} \mathcal{G}}  = \frac{1}{2} \frac{\partial}{\partial x^{i}} \wedge \frac{\partial}{\partial x^{i}} - \frac{1}{2} C_{\alpha \beta} ^{\gamma} \eta_{\gamma} \frac{\partial}{\partial \eta_{\alpha}} \wedge \frac{\partial}{\partial \eta_{\beta}}. \nonumber
 \end{equation}
On the other hand, let $h: T_{\xi^\prime}^{\ast} O(\xi^{\prime}) \longrightarrow \mathcal{R}$ be the Hamiltonian function on $T_{\xi^\prime}^{\ast} O(\xi^{\prime}).$ So, its Hamiltonian equations according to equation (\ref{HEh}) are as follows:
\begin{equation}
\label{HE2}
\frac{d \eta_{\alpha}}{dt} = -  C_{\alpha \beta} ^{\gamma} \eta_{\gamma}\frac{\partial h}{\partial \eta_{\beta}}.
\end{equation}
Suppose $\lbrace~,~ \rbrace_{K.K} $  is the symbol of Kirillov-Kostant bracket on $C^{\infty} ( T^{\ast} _{\xi^{\prime}} O(\xi^{\prime} )).$  By (\ref{kk}) we have
\begin{equation}
\label{hvf}
 \mathcal{X}_{h} (f ) = \Pi (d f , d h ) = \lbrace f, h \rbrace_{K.K}.
\end{equation}
 Furthermore, as we proved in \cite{HA1}, for  Hamiltonian $H=(p, h): M \times T_{\xi^{\prime}} ^{\ast} O(\xi^{\prime})  \longrightarrow \mathbb{R},$ we have
 \begin{equation}
 \label{bra1}
 \lbrace F, H \rbrace _{A^{\ast} \mathcal{G} } = \lbrace f, h \rbrace_{K.K},
 \end{equation}
 where $ F = (p , f ) \in C^{\infty} (M \times T^{\ast} _{\xi^{\prime}} O(\xi^{\prime} )) $, $ f \in C^{\infty} (T^{\ast} _{\xi^{\prime}} O(\xi^{\prime} ))$ and $p\in M.$ Therefore, according to equations (\ref{HVF}), (\ref{hvf}) and (\ref{bra1}), we conclude that
  \begin{equation}
 \label{equality}
  \mathcal{X}_{H} ^{\Pi_{A^{\ast} \mathcal{G}}} (F) =\mathcal{X}_{h} (f ),
 \end{equation}
where $ \mathcal{X}_{H} ^{\Pi_{A^{\ast} \mathcal{G}}}$ is Hamiltonian vector field associated to Hamiltonian $H: M \times T_{\xi^{\prime}} ^{\ast} O(\xi^{\prime})  \longrightarrow \mathbb{R}.$ 
\\

If we consider the Hamiltonian function $H: M \times T_{\xi^{\prime}} ^{\ast} O(\xi^{\prime})  \longrightarrow \mathbb{R}, H = (p, h),$ and $\delta = (p,\lambda) \in M \times T_{\xi^{\prime}} ^{\ast} O(\xi^{\prime}),$ we have that $H (p, \lambda ) = h(\lambda),$
where $h: T_{\xi^\prime}^{\ast} O(\xi^{\prime}) \longrightarrow \mathbb{R}$ is Hamiltonian function on $T_{\xi^\prime}^{\ast} O(\xi^{\prime}).$\\

 Now let we have a control system $\sigma: C \to A\mathcal{G}=M \times T_{\xi^\prime}O(\xi^\prime)$, $\sigma(c)=\left(p, \sigma_1(c)\right)$ for all $c\in C$, where $\pi: C\to M$ is control space and $\sigma_1:C \to T_{\xi^\prime}O(\xi^\prime)$. Let us consider the Hamiltonian $H\in C^\infty \left(A^\ast\mathcal{G}\times _M C\right)$ as follows: 
 $$H(\eta, c)=\langle\eta, \sigma(c)\rangle-L(c),$$
 where $L:C\to \mathbb{R}$ is the cost function. For every $\eta=(p, \eta_1) \in A^\ast\mathcal{G}=M \times T^\ast_{\xi^\prime}O(\xi^\prime)$,  
 where $\eta_1 \in T^\ast_{\xi^\prime}O(\xi^\prime),$ we have $H(\eta, c)=\langle\eta_1, \sigma_1(c)\rangle-L(c).$
  So $H=(p, h)$, where $h\in C^\infty \left(T_{\xi^\prime}O(\xi^\prime)\times  C\right)$ and $h(\eta_1, c) =\langle\eta_1, \sigma_1(c) \rangle-L(c).$
 
As seen above, equations (\ref{HE1}), the Hamiltonian equations for $H$ is as follows: 
 \begin{equation}
  \frac{d x^{i}}{dt} = \frac{\partial H}{\partial\eta_\alpha},  \quad \frac{d \eta_{\alpha}}{dt} = -\left(\frac{\partial H}{\partial x^{i}}+  C_{\alpha \beta}  ^{\gamma} \eta_{\gamma}\frac{\partial H}{\partial \eta_{\beta}}\right).\nonumber
 \end{equation}
 So by using the  equation (\ref{critical tra}) we obtain the equations for the critical trajectories as

\begin{eqnarray}
0&=&\frac{\partial H}{\partial\eta_\alpha},\nonumber\\
\dot{\eta_\alpha}&=&- \left(\frac{\partial H}{\partial x_i}+ \eta_\gamma C^\gamma_{\alpha\beta}\frac{\partial H}{\partial \eta_\beta}\right),\nonumber\\
0&=&\frac{\partial H}{\partial u^c}.
\label{compl1}
\end{eqnarray}

For a Hamiltonian function  $H=(p, h)$ on $M \times T_{\xi^{\prime}} ^{\ast} O(\xi^{\prime}),$ according to the  equations (\ref{HE1}), (\ref{HE2}) and (\ref{equality}), the equations for critical trajectories will be 
\begin{eqnarray}
\dot{\eta_\alpha}&=&-  \eta_\gamma C^\gamma_{\alpha\beta}\frac{\partial h}{\partial \eta_\beta}, \nonumber\\
0&=&\frac{\partial h}{\partial u^c}.
\label{simp1}
\end{eqnarray}
\begin{conclusion}
The results of this work can be  significant in the control theory because our reduction in the new Lie groupoids as well in
 Lie algebroids, i.e. co-adjoint Lie algebroid significantly simplify the Hamiltonian equations associated with the control system. As we see, in the illustrated example, by using the reduction in the co-adjoint Lie algebroid, one can easily
reach the optimal solutions of the system. In other words, in the case of trivial groupoid, instead of finding the solutions of the more complicated Hamiltonian system (\ref{compl1}) one can consider the solutions of the simple Hamiltonian system (\ref{simp1}). To do so, we consider
the trivial Lie groupoid and an optimal control problem on its co-adjoint Lie algebroid and show that the
optimal control problem can be reduced to the optimal control problem on the co-tangent bundle of the orbits
of the co-adjoint representation of the Lie group. Also, we show that the extermal solutions of the optimal
control problem on the co-adjoint Lie algebroid of the trivial Lie groupoid are obtained from the solutions of
the corresponding Hamiltonian system on the co-tangent bundle of the co-adjoint orbits of the Lie group.  
\end{conclusion}
{\bf Statements and Declarations:}\\
{\bf Financial conflict:}\\
There is no any financial and non-financial conflict of interest associated with this
work.\\
{\bf Data Availability Statement:}\\
Data sharing not applicable to this article as no datasets were generated or analysed during the current study.\\
{\bf  Funding:}\\
The author did not receive support from any organization for the submitted work.\\
{\bf Competing Interests:}\\
There is no any  competing interest for submitted work.

\small{

}


\begin{thebibliography}{99}
\bibitem{HA1} Gh. Haghighatdoost and R. Ayoubi, \textit{Hamiltonian systems on co-adjoint Lie groupoids}, Journal of Lie Theory 31, 2 (2021), 493-516. 
\bibitem{HA2} Gh. Haghighatdoost and R. Ayoubi, \textit{Euler integrable Hamiltonian system on co-adjoint Lie groupoids}, researchGate DOI: 10.13140/RG.2.2.25372.62081.(2021)
\bibitem{HA3} Gh. Haghighatdoost and R. Ayoubi, \textit{Generalized geometric Hamilton-Jacobi theorem on Lie algebroids}, arXiv:1902.06969v1, (2019).
\bibitem{Joz} Michal Jozwikowski, \textit{Optimal control theory on almost Lie algebroids}, PhD dissertation, Institute of Mathematics
Polish Academy of Sciences,arXiv:1111.1549v1 [math.OC] 7 Nov 2011.
\bibitem{GJ}  Janusz Grabowski, Micha l Jozwikowski,  \textit{Pontryagin Maximum Principle, a generalization},January 2011 SIAM Journal on Control and Optimization 49(3):1306-1357 
DOI: 10.1137/090760246 and arXiv:0905.2767v2 [math.OC] 1 Dec 2011.
\bibitem{Mac1} K.C.H. Mackenzie, \textit{General theory of Lie groupoids and Lie algebroids},
 London Math. Soc. Lecture notes series 213, Cambridge University Press, Cambridge, (2005).
\bibitem{Mar} Eduardo Martinez \textit{Reduction in optimal control theory}, Reports on Mathematical Physics, No.1, Vol.53, (2004).
\bibitem{BW} J. Baillieul, J.C. Willems Editors \textit{Mathematical Control Theory}, 1999 Springer Science+Business Media New York
Originally published by Springer-Verlag New York, Inc. in 1999.
\bibitem{suss} H. J. Sussmann, \textit{Geometry and Optimal Control}, Mathematical Control Theory" (a book of essays a in honor
of Roger W. Brockett on the o ccasion of his 60th Birthday), J. Baillieul and J. C. Willems
Eds., Springer-Verlag, 1998 (ISBN 0-387-98317-1), pages 140-198.
\bibitem{Jur} V. Jurdjevic, \textit{Geometric control theory}, Cambridge University Press 1997.
\bibitem{Jur2} V. Jurdjevic: \textit{Optimal control problems on Lie groups}: Crossroads between geometry and mechanics, In
Geometry of Feedback and Optimal Control, B. Jakubczyk and W. Respondek (eds.), Marcel-Dekker,
1993.
\bibitem{DeL} M. de Le´on, J. C. Marrero, E. Mart´ınez: \textit{Lagrangian submanifolds and dynamics
on Lie algebroids}, J. Phys. A : Math. Gen, 38 (2005), 241–308.
\bibitem{Marr} J. C. Marrero: \textit{Hamiltonian dynamics on Lie algebroids, unimodularity and
preservation of volumes}, arXiv preprint arXiv:0905.0123v1, (2009).
\bibitem{Bos}
R. Bos:
\newblock {\em Geometric quantization of Hamiltonian actions of Lie algebroids and Lie groupoids},
\newblock Int. J. Geom. Methods Mod. Phys, {\bf 4} (2007), 389--436.
\end{thebibliography}
\end{document}